\documentclass[reqno]{article}

\usepackage{xcolor}
\usepackage{amssymb}
\usepackage{amsfonts}
\usepackage{amsmath}
\numberwithin{equation}{section}
\usepackage{enumerate}
\usepackage[pagewise]{lineno}

\newtheorem{theorem}{Theorem}[section]
\newtheorem{proposition}[theorem]{Proposition}
\newtheorem{lemma}[theorem]{Lemma}
\newtheorem{corollary}[theorem]{Corollary}
\newtheorem{definition}[theorem]{Definition}
\newtheorem{example}[theorem]{Example}

\newtheorem{remark}[theorem]{Remark}

\title {Reducibility of self-adjoint linear relations and application to generalized Nevanlinna functions}
 \author{Muhamed Borogovac, \\
 muhamed.borogovac@gmail.com}
\date{September 14, 2020} 
\begin{document}
\maketitle

\textbf{Abstract}  

Necessary and sufficient conditions for reducidibility of a self-adjoint linear relation in a Krein space are given. Then a generalized Nevanlinna function $Q$, represented by a self-adjoint linear relation $A$, is decomposed by means of the reducing subspaces of $A$. The sum of two functions $Q_{i}{\in N}_{\kappa_{i}}\left( \mathcal{H} \right),\thinspace i=1,\thinspace 2$, minimally represented by the triplets $\left( \mathcal{K}_{i},A_{i},\Gamma_{i} \right)$, is also studied. For that purpose, a model $( \tilde{\mathcal{K}},\tilde{A},\tilde{\Gamma } )$ to represent $Q:=Q_{1}+Q_{2}$ in terms of $\left( \mathcal{K}_{i},A_{i},\Gamma_{i} \right)$ is created. By means of that model, necessary and sufficient conditions for $\kappa =\kappa_{1}+\kappa_{2}$ are proven in analytic terms. At the end, it is explained how degenerate Jordan chains of the representing relation $A$ affect reducing subspaces of $A$ and decomposition of the corresponding function $Q$. 



\section{Preliminaries and introduction}\label{s2}

\textbf{1.1} \textbf{Preliminaries} Let $N$, $R$, and $C$ denote sets of positive integers, real numbers, and complex numbers, respectively. Let $(.,.)$ denote the (definite) scalar product in the Hilbert space $\mathcal{H}$, and let $\mathcal{L}(\mathcal{H})$ denote the space of bounded linear operators in $\mathcal{H}$.

\begin{definition}\label{definition2} An operator-valued complex function $Q:\mathcal{D}\left(Q \right)\to \mathcal{L}(\mathcal{H})$ belongs to the class of generalized Nevanlinna functions $N_{\kappa}\left(\mathcal{H} \right)$ if it satisfies the following requirements:

\begin{itemize}
\item $Q$ is meromorphic in $C\thinspace \backslash \thinspace R$,
\item ${Q\left(\bar{z} \right)}^{\ast}=Q\left(z \right),\thinspace z\in \mathcal{D}\left(Q \right)$, and
\item the Nevanlinna kernel $N_{Q}\left(z,w \right):=\frac{Q\left(z \right)-{Q\left(w \right)}^{\ast}}{z-\bar{w}};\thinspace z,\thinspace w \in \mathcal{D}(Q)\cap C^{+}$,
\end{itemize}
has $\kappa$ negative squares. In other words, for arbitrary $n\in N, z_{1}, \mathellipsis, z_{n}\in \mathcal{D}(Q)\cap C^{+}$ and $h_{1}, \mathellipsis,h_{n}\in \mathcal{H}$, the Hermitian matrix $\left(N_{Q}\left(z_{i},z_{j} \right)h_{i},h_{j} \right)_{i,j=1}^{n}$ has $\kappa$ negative eigenvalues at most, and for at least one choice of $n; z_{1}, \mathellipsis, z_{n}$, and $h_{1}, \mathellipsis,h_{n}$, it has exactly $\kappa $ negative eigenvalues.
\end{definition}
It is easy to verify that Nevanlinna kernel is a \textit{Hermtian} kernel, i.e. $N_{Q}\left(z,w \right)^{*}=N_{Q}\left(w,z \right), \thinspace z,\thinspace w \in \mathcal{D}(Q)\cap C^{+}$. 

The following definitions of a linear relation and basic concepts related to it can be found for example in \cite{A,BHS,S}. In the sequel, $\mathcal{H}$, $\mathcal{K}$, $\mathcal{M}$ are inner product spaces. 
Recall, a set $M$ is called \textit{linear manifold} (or \textit{linear space}) if for any two vectors $ x, y \in M $ and for any two scalars $\alpha, \beta \in C $ it holds $ \alpha x + \beta y \in M $. A \textit{linear relation} from $\mathcal{H}$ into $\mathcal{K}$ is a linear manifold $T$ of the product space $\mathcal{H}\times \mathcal{K}$. If $\mathcal{H}=\mathcal{K}$, $T$ is said to be a \textit{linear relation in} $\mathcal{K}$. We will use the following concepts and notations for linear relations, $T$ and $S$ from $\mathcal{H}$ into $\mathcal{K}$ and a linear relation $R$ from $\mathcal{K}$ into $\mathcal{M}$. 
\[
D\left( T \right):=\left\{ f\in \mathcal{H}\vert \left\{ f,g \right\}\in T\, for\, some\, g\in \mathcal{K} \right\},
\]
\[
R\left( T \right):=\left\{ g\in \mathcal{K}\vert \left\{ f,g \right\}\in T\, for\, some\, f\in \mathcal{H} \right\},
\]
\[
\ker T:=\left\{ f\in \mathcal{H}\vert \left\{ f,0 \right\}\in T \right\},
\]
\[
T(0):=\left\{ g\in \mathcal{K}\vert \left\{ 0,g \right\}\in T \right\},
\]
\[
T\left( f \right):=\left\{ g\in \mathcal{K}\vert \left\{ f,g \right\}\in T \right\},\, \, (f\in D\left( T \right)), 
\]
\[
T^{-1}:=\left\{ \left\{ g,f \right\}\in \mathcal{K} \times \mathcal{H}\vert \left\{ f,g \right\}\in T \right\},
\]
\[
zT:=\left\{ \left\{ f,zg \right\}\in \mathcal{H} \times \mathcal{K}\vert \left\{ f,g \right\}\in T \right\},\, \, (z\in C),
\]
\[
S+T:=\left\{ \{f,g+k\}\vert \left\{ f,g \right\}\in S,\left\{ f,k \right\}\in T \right\},
\]
\[
RT:=\left\{ \{f,k\}\in \mathcal{H} \times \mathcal{M}\vert \left\{ f,g \right\}\in T,\left\{ g,k \right\}\in R\, for\, some\, g\in \mathcal{K} \right\},
\]
\[
T^{+}:=\left\{ \{k,h\}\in \mathcal{K} \times \mathcal{H} \vert \left[ k,g \right]=\left( h,f \right)\, for\, all\, \left\{ f,g \right\}\in T \right\},
\]
\[
T_{\infty }:=\left\{ \left\{ 0,g \right\}\in T \right\}.
\]
Note that in definition of the adjoint linear relation $T^{+}$, we use the following notation for inner product spaces: $\left( \mathcal{H}, (.,.) \right)$ and $\left( \mathcal{K}, [.,.] \right)$. 

If $mul T:=T(0)=\lbrace 0\rbrace$, we say that $T$ is an \textit{operator}, or \textit{single-valued} linear relation. A linear relation is \textit{closed} if it is a closed subset in the product space $\mathcal{H} \times \mathcal{K}$. 

Let $A$ be a linear relation in $\mathcal{K}$. We say that $A$ is \textit{symmetric} (\textit{self-adjoint}) if it holds $A\subseteq A^{+}$ ($A=A^{+})$. Every point $\alpha \in C$ for which $\left\{ f,\alpha f \right\}\in A$, with some $f\ne 0$, is called a \textit{finite eigenvalue}. The corresponding vectors are \textit{eigenvectors} belonging to the eigenvalue $\alpha $. The set that consists of all points $z\in C$ for which the relation $\left( A-zI \right)^{-1}$ is an operator defined on the entire $\mathcal{K}$, is called the \textit{resolvent} set $\rho (A)$.
\\

Let $\kappa \in N\cup\lbrace 0\rbrace$ and let $\left(\mathcal{K},\thinspace \left[
.,. \right] \right)$ denote a \textit{Krein space}. That is a complex vector space on which a scalar product, i.e. a Hermitian sesquilinear form $\left[ .,. \right]$, is defined such that the following decomposition 
\[
\mathcal{K}=\mathcal{K}_{+} \dot{+} \mathcal{K}_{-}
\]
of $\mathcal{K}$ exists, where $\left(\mathcal{K}_{+},\left[ .,. \right] \right)$ and $\left(\mathcal{K}_{-},-\left[ .,. \right] \right)$ are Hilbert spaces which are mutually orthogonal with respect to the form $\left[ .,. \right]$. Every Krein space $\left( \mathcal{K},\, \left[ .,.\right] \right)$ is \textit{associated} with a Hilbert space $\left( \mathcal{K},\, \left( .,.\right) \right)$, which is defined as a direct and orthogonal sum of the Hilbert spaces $\left(\mathcal{K}_{+},\left[ .,. \right] \right)$ and $\left(\mathcal{K}_{-},-\left[ .,. \right] \right)$. The topology in a Krein space $\mathcal{K}$ is the topology of the associated Hilbert space $\left( \mathcal{K},\, \left( .,.\right) \right)$. For properties of Krein spaces see e.g. \cite[Chapter V]{Bog}. 

If the scalar product $\left[.,. \right]$ has $\kappa \thinspace (<\infty)$ negative squares, then we call it a
\textit{Pontryagin space} of the index $ \kappa $. The definition of a Pontryagin space and other concepts related to it can be found e.g. in \cite{IKL}. 
\\

The following construction of a Pontryagin space can be found in \cite{HSW,KL1,DLS} and a similar construction of a Hilbert space can be found in \cite{LT}:

For any generalized Nevanlinna function $Q$, a linear space $L(Q)$ with a (possibly degenerate) indefinite inner product $\left[.,. \right]$ can be introduced as follows:

Consider the set of all \textbf{finite} formal sums
\[
\sum {\varepsilon_{z}h_{z}}, z\in \mathcal{D}\left(Q \right),
\]
where $h_{z}\in \mathcal{H}$, and $\varepsilon_{z}$ is a symbol associated with each $z\in \mathcal{D}\left(Q \right)$. Then, an inner product is defined by:
\[
\left[\varepsilon_{z}h_{z},\varepsilon_{w}h_{w} \right]:=\left(\frac{Q\left(z \right)-Q\left(\bar{w} \right)}{z-\bar{w}}h_{z},h_{w} \right), z, w\in \mathcal{D}\left(Q \right), \thinspace \thinspace z\ne
\bar{w},\thinspace h_{z}, h_{w}\in \mathcal{H},
\]
\[
\left[\varepsilon_{z}h_{z},\varepsilon_{\bar{z}}h_{\bar{z}} \right]:=\left(Q^{'}\left(z \right)h_{z},h_{\bar{z}} \right), z\in \mathcal{D}\left(Q \right).
\]

To ease communication, let us call $L(Q)$ the \textit{state manifold of} $Q$.
The linear relation defined by
\[
A_{0}:=l.s.\left\{\left\{\sum\limits_s {\varepsilon_{z_{s}}h_{s}},\sum\limits_s {z_{s}\varepsilon_{z_{s}}h_{s}} \right\}: \sum\limits_s h_{s} =0, z_{s}\in \mathcal{D}\left(Q \right)\right\}
\]
is symmetric. For $z_{0} \in \mathcal{D} \left(Q \right)$, the operator $ \Gamma_{z_{0}}: \mathcal{H}\rightarrow L \left(Q \right)$ is defined by $\Gamma_{z_{0}}h= \varepsilon_{z_{0}}h$.
The Pontryagin space $\mathcal{K}$ is obtained by factorization of $L(Q)$ with respect to its isotropic part $L_{0}^{0}:=L(Q)\cap {L(Q)\mathrm{\thinspace}}^{[\bot]}$ and by completion of the factor space. It is called the \textit{state space} of $Q$. In the process, $A_{0}$ and $\Gamma_{z_{0}}$ give rise to the self-adjoint relation $A$ in $\mathcal{K}$ and bounded linear operator $\Gamma :\mathcal{H}\rightarrow \mathcal{K}$, with $z_{0} \in \rho(A)$. Then the following theorem holds.

\begin{theorem}\label{theorem1} A function $ Q:\mathcal{D}(Q)\to L(\mathcal{H})$ is a generalized Nevanlinna function of the index $\kappa $, denoted by $Q\in N_{\kappa}(\mathcal{H})$, if and only if it has a representation of the form

\begin{equation}
\label{eq2}
Q\left(z \right)={Q(z_{0})}^{\ast}+(z-\bar{z_{0}})\Gamma^{+}\left(I+\left(z-z_{0} \right)\left(A-z \right)^{-1} \right)\Gamma, z\in \mathcal{D}\left(Q \right),
\end{equation}
where $A$ is a self-adjoint linear relation in some Pontryagin space $(\mathcal{K}, [.,.])$ of the index $\tilde{\kappa}\ge \kappa; \Gamma:\mathcal{H}\to \mathcal{K}$ is a bounded operator. (Obviously $\rho (A) \subseteq \mathcal{D}(Q)$.) This representation can be chosen to be minimal, that is,
\begin{equation}
\label{eq2a}
\mathcal{K}=c.l.s.\left\{\Gamma_{z}h:z\in \rho \left(A \right),h\in \mathcal{H} \right\},
\end{equation}
where
\begin{equation}
\label{eq2b}
\Gamma_{z}:=\left(I+\left(z-z_{0} \right)\left(A-z \right)^{-1} \right)\Gamma .
\end{equation}

If realization (\ref{eq2}) is minimal, then $Q\in N_{\kappa}(\mathcal{H})$ if and only if $\tilde{\kappa}$ \textit{equals} $\kappa $. In the case of minimal representation $\rho(A)=D(A)$ and the triple $(\mathcal{K},\thinspace A,\thinspace \Gamma)$ is uniquely determined (up to isomorphism).
\end{theorem}

Such operator representations were developed by M. G. Krein and H. Langer \cite{KL1,KL2} and later converted to representations in terms of linear relations (see e.g. \cite{DLS, HSW}).

In this paper, a point $\alpha \in C$ is called a \textit{generalized pole} of $Q$ if it is an eigenvalue of the representing relation $A$. It may be an isolated singularity, i.e. an ordinary pole, as well as an embedded singularity of $Q$. The latter may be the case only if $\alpha \in R$.
\\

\textbf{1.2. Introduction} We start Section \ref{s4} with extending the definition of reducibility of operators in Hilbert spaces to reducibility of linear relations in Krein spaces. Then in Lemma \ref{lemma6} we prove several statements about decompositions, i.e. about relation matrix, of a linear relation in a Krein space $\mathcal{K}$ that we need in the proof of the main result, Theorem \ref{theorem2}. In that theorem we give necessary and sufficient conditions for a self-adjoint linear relation $A$ in $\mathcal{K}$ to be reduced to the sum $A=A_{1}\left[ + \right]A_{2}$, where ``$\left[ + \right]$'' is direct and orthogonal sum of linear relations, $A_{i}$ are self-adjoint linear relations in the reducing subspaces $\mathcal{K}_{i}$, and $\mathcal{K}=\mathcal{K}_{1}\left[ \dotplus \right]\mathcal{K}_{2}$. Then, by means of reducing subspaces and reducing linear relations we study decompositions of a generalized Nevanlinna function $Q$. 

The number of negative squares $\kappa \in N \cup \lbrace 0 \rbrace $ is an important feature of the generalized Nevanlinna function $Q$. Recall that, if functions $Q_{i},\thinspace i=1,\thinspace 2$, satisfy 

\begin{enumerate}[(i)]
\item $Q_{i}{\in N}_{\kappa_{i}}(\mathcal{H})\thinspace ,\thinspace 0\le \kappa_{i}, i=1,\thinspace 2,$ 
\item $Q\left( z \right)=Q_{1}\left( z \right)+Q_{2}\left( z \right),$
\end{enumerate}

then $Q$ belongs to some generalized Nevanlinna class $N_{\kappa }\left( \mathcal{H} \right)$ and $ \kappa \le \kappa_{1}+\kappa_{2}$ holds. 

There are two basic questions: 

\begin{enumerate}[(a)]
\item Given function $Q{\in N}_{\kappa }\left( \mathcal{H} \right)$, under what conditions does there exist a decomposition \\
$Q\left( z \right)=Q_{1}\left( z \right)+Q_{2}\left( z \right),\thinspace Q_{i}{\in N}_{\kappa_{i}}\left( \mathcal{H} \right),\thinspace i=1,\thinspace 2$, that satisfies $\kappa =\kappa_{1}+\kappa_{2}$? 
\item Given two functions $Q_{i}{\in N}_{\kappa_{i}}(\mathcal{H})\thinspace ,\thinspace i=1,\thinspace 2$, is the number of negative squares preserved in the sum $Q=Q_{1}+Q_{2}$ or not? 
\end{enumerate}

In other words, we investigate the circumstances under which functions 
${Q,\thinspace Q}_{1}$ and $Q_{2}$ that satisfy (i) and (ii) also satisfy

\begin{enumerate}[(iii)]
\item $\kappa_{1}+\kappa_{2}=\kappa$. 
\end{enumerate}

The question of preservation of the number of negative squares of the sum of Hermitian kernels $K(z,w)=K_{1}(z,w)+K_{2}(z,w)$ was studied in \cite{ADRS}. The authors give necessary and sufficient conditions for $\kappa_{1}+\kappa_{2}=\kappa$ in terms of complementary reproducing kernel Pontryagin spaces $\mathcal{K}_{1}, \mathcal{K}_{2}$, c.f. \cite[Theorem 1.5.5]{ADRS}. We alternatively give necessary and sufficient conditions for $\kappa_{1}+\kappa_{2}=\kappa$ in terms of triplets $\left( \mathcal{K}_{i},A_{i},\Gamma_{i} \right), \thinspace i=1, 2$, associated with minimal representations of the form (\ref{eq2}), c.f. Theorem \ref{theorem6}. 

The question of preservation of the number of negative squares in products, sums, and in some transformations of generalized Nevanlina functions has been, among other topics, summarised in the survey \cite{Lu}. In the present paper we prove analytic criteria that establish whether the sum of the indexes of the functions that comprise the sum is equal or it is greater than the negative index of the sum.  

It is very difficult to determine the negative index $\kappa$ of a given generalized Nevanlinna function. The established relation between negative indexes of the above sum (ii) gives us information that might help in determining the numbers of negative indexes of the functions in the sum. 

There are interesting results about decompositions of generalized Nevanlinna 
functions in \cite{DL, KL2}, for matrix and scalar functions represented by unitary 
and self-adjont operators. In those papers, the decompositions of $Q\thinspace \in 
\thinspace N_{\kappa }^{n \times n}$ were obtained by means of spectral families of 
the representing operators and their appropriate invariant spectral 
subspaces. The decomposing functions $Q_{i}$ obtained by that method must 
have disjoint sets of generalized poles (see \cite[Proposition 3.1]{DL}). In the 
present article, we do not use spectral families and spectral 
subspaces; we use instead a concept of the reducing subspaces of the representing self-adjoint 
relation in the Pontryagin state space. That way we obtain decomposition where decomposing functions $Q_{i},\thinspace i=1,\thinspace 2$ may have common generalized poles. 

In Theorem \ref{theorem3}, we give a general answer on the question (a); we 
decompose function $Q$ by means of reducing subspaces $\mathcal{K}_{i}$ and reducing 
relations $A_{i}$ of the representing relation $A$. 

Regarding sums of generalized Nevanlina functions, in \cite[Proposition 3.2]{DL} it 
has been proven that the sum of two generalized Nevanlinna matrix functions 
preserves the number of negative squares under the condition that 
functions in the sum have disjoint sets of generalized poles. In our study we do not use that condition. 

We start the study of the sum $Q:=Q_{1}+Q_{2}$ with two functions $Q_{i}{\in 
N}_{\kappa_{i}}\left( \mathcal{H} \right),\thinspace i=1,\thinspace 2,$ represented 
minimally in Pontryagin spaces by triplets $\left( \mathcal{K}_{i},A_{i},\Gamma_{i} 
\right),\thinspace \thinspace i=1,\thinspace 2$. Then we create Pontryagin 
space $\tilde{\mathcal{K}}:=\mathcal{K}_{1}\left[\dotplus \right]\mathcal{K}_{2}$, and representation of the 
function $Q:=Q_{1}+Q_{2}$ in terms of the triplets $\left( 
\mathcal{K}_{i},A_{i},\Gamma_{i} \right)$. That representation, denoted by (\ref{eq36}) in 
the text, we call \textit{orthogonal sum representation.} Then, in Theorem \ref{theorem6}, we describe the structure of the possibly non-minimal state space $\tilde{\mathcal{K}}:=\mathcal{K}_{1}\left[\dotplus 
\right]\mathcal{K}_{2}$ representing the sum $Q=Q_{1}+Q_{2}$. In Corollary \ref{corollary5}, we 
give necessary and sufficient conditions for $\kappa \thinspace =\kappa 
_{1}+\kappa_{2}$ in terms of the inner structure of the state space 
$\tilde{\mathcal{K}}$.

In Theorems \ref{theorem8} and \ref{theorem10} we prove some analytic criteria for $\kappa \thinspace =\kappa_{1}+\kappa_{2}$ or $\kappa <\kappa_{1}+\kappa _{2}$. These criteria are easy to use; we do not need to know operator representations of the functions comprising the sum. Given how Definition 
\ref{definition2} is impractical for use and how difficult it is to find operator 
representations, our criteria are useful tool for research of both, the 
underlying state space, and features of the sum $Q\thinspace :=\thinspace Q_{1}+Q_{2}$. 

In Proposition \ref{proposition6}, we decompose a function $Q$ by means of Theorem \ref{theorem3} using linear spans of non-degenerate Jordan chains as reducing subspaces. Proposition \ref{proposition6} is a straightforward result that we needed to approach the more complicated case of degenerate chains 
which we study in Proposition \ref{proposition8}. In Proposition \ref{proposition8} we consider the model where the self-adjoint operator $A$ in a Pontryagin space $\mathcal{K}$ has two simple, 
independent, and degenerate chains (neutral eigenvectors) at $\alpha 
\thinspace \in \thinspace R$. We prove that, unlike non-degenerate chains, 
studied in Proposition \ref{proposition6}, the two degenerate chains at $\alpha 
\thinspace \in \thinspace R$ cannot reduce the representing operator and 
cannot induce two different functions $Q_{i}$ in any decomposition of $Q$. The conclusion of Section \ref{s10} is in Corollary \ref{corollary10}

\section{Reducing subspaces of the self-adjoint linear relation in the Krein space}\label{s4}
In the sequel ``$\left[ + \right]$'', rather than "$\left[ \dotplus \right]$", denotes direct and orthogonal sum of both, relations and vectors. From the context it is usually clear when we deal with "operator-like" addition of linear relations, as well as we deal with addition of relations as subspaces, and addition of vectors. If necessary, we will specify. 

\begin{lemma}\label{lemma4} Assume that $\mathcal{K}_{1}$ and $\mathcal{K}_{2}$ are Krein spaces and 
$A_{l}\subseteq \mathcal{K}_{l}^{2},\thinspace l=1,\thinspace 2$, are linear 
relations. We can define direct orthogonal sum
\[
\mathcal{K}:=\mathcal{K}_{1}\left[ + \right]\mathcal{K}_{2}
\]
and
\[
A=A_{1}\left[ + \right]A_{2}:=\left\{ \left( {\begin{array}{*{20}c}
h_{i} \left[ + \right] h_{j}\\
h_{i}^{'} \left[ + \right] h_{j}^{'}\\
\end{array} } \right ) : \left( {\begin{array}{*{20}c}
h_{l}\\
h_{l}^{'}\\
\end{array} } \right)\in A_{l},\thinspace \thinspace l=1,\thinspace 2 
\right\}\subseteq \mathcal{K}^{2}.
\]
The linear relation $A:=A_{1}\left[ + \right]A_{2}$ is symmetric (self-adjoint) in $\mathcal{K}$ if and only if linear relations $A_{l}\subseteq \mathcal{K}_{l}^{2},\thinspace \thinspace l=1,\thinspace 2$, are symmetric (self-adjoint).
\end{lemma}
\textbf{Proof:} This lemma is a straightforward verification and left to the reader. $\square$ 

Let $A$ be a linear relation in the Krein space $\mathcal{K}$, 
\[
\mathcal{K}:=\mathcal{K}_{1}\left[ + \right]\mathcal{K}_{2}
\]
where nontrivial subspaces $\mathcal{K}_{l}$ are also Krien spaces and $E_{l}:\mathcal{K}\to \mathcal{K}_{l},\thinspace l=1,\thinspace 2$, are the corresponding orthogonal projections. The following four linear relations can be introduced 
\[
A_{i}^{j}:=\left\{ \left( {\begin{array}{*{20}c}
h_{i}\\
h_{i}^{j}\\
\end{array} } \right):h_{i}\in D\left( A \right)\cap \mathcal{K}_{i},h_{i}^{j}\in 
E_{j}A\left( h_{i} \right) \right\}\subseteq \mathcal{K}_{i} \times \mathcal{K}_{j},\thinspace 
i,j=1,\thinspace 2.
\]
In this notation the subscript ``$i$'' is associated with the domain 
subspace $\mathcal{K}_{i}$, the superscript ``$j$'' is associated with the range 
subspace $\mathcal{K}_{j}$. For example $\left( {\begin{array}{*{20}c}
h_{1}\\
h_{1}^{2}\\
\end{array} } \right)\in A_{1}^{2}$.

Let us now extend the definition of the reducing subspaces of the unbounded 
operator in the Hilbert space, see e.g. \cite[Section 40]{AG}, to the reducing 
subspaces of the (multivalued) linear relation in Krein space. 

\begin{definition}\label{definition4} Let $\left( \mathcal{K},\left[ .,. \right] \right)$ be a Krein 
space and let $\mathcal{K}_{1}\subset \mathcal{K}$ be a nontrivial Krein subspace of $\mathcal{K}$, and $\mathcal{K}_{2}=\mathcal{K}\left[ - \right]\mathcal{K}_{1}$. We will say that the subspaces $\mathcal{K}_{1}$ and $\mathcal{K}_{2}$ reduce relation $A$ if there exist linear relations $A_{i} \subseteq \mathcal{K}_{i} \times \mathcal{K}_{i},\thinspace i=1,\thinspace 2$, such that it holds 
\[
A=A_{1}\left[ + \right]A_{2},
\]
where $[+]$ stands for direct orthogonal addition of relations, as defined in Lemma \ref{lemma4}. The relations $A_{i}$ are called reducing relations of $A$.
\end{definition}

Recall, if $\mathcal{K}$ is a Pontryagin space and $\mathcal{K}_{1}$ is a non-degenerate closed subspace, then $\mathcal{K}=\mathcal{K}_{1}[+]\mathcal{K}_{2}$ and both $\mathcal{K}_{i}, \thinspace i=1, 2$, are also Pontryagin spaces, see \cite[Theorem 3.2 and Corollary 2]{IKL}. 

\begin{lemma}\label{lemma6} Let $A$, $A_{i}^{j},\thinspace \thinspace \mathcal{K}_{i},\thinspace \thinspace E_{i};\thinspace i,\thinspace j=1,\thinspace 2$ be introduced as above. If for either of orthogonal projections $E_{i}:\mathcal{K} \to \mathcal{K}_{i},\thinspace i=1,\thinspace 2,$ it holds $E_{i}(D(A))\subseteq D(A)$, then 

\begin{enumerate}[(i)]
\item $E_{1}\left( A\left( 0 \right) \right)=A_{1}^{1}\left( 0 \right)=A_{2}^{1}\left( 0 \right),\thinspace \thinspace E_{2}\left( A\left( 0 \right) \right)=A_{1}^{2}\left( 0 \right)=A_{2}^{2}\left( 0 \right)$.

\item $A=\left( A_{1}^{1}+A_{1}^{2} \right)\hat{+}\left( A_{2}^{1}+A_{2}^{2} \right),$ where ``$+$'' stands for operator-like addition, and ``$\hat{+}$'' stands for addition of the subspaces, not necessarily direct. 

\item If $B: \mathcal{K}_{i}\rightarrow \mathcal{K}_{j}, \thinspace i,j = 1, 2$, is a closed relation, then 
\begin{equation}
\label{eq22}
B \left( 0 \right)=D\left(B^{[\ast]}\right)^{[\bot ]} \left( \subseteq \mathcal{K}_{j}\right).
\end{equation}
\item If $A$ is symmetric, then it holds: $A_{2}^{1}\subseteq {A_{1}^{2}}^{\left[ \ast \right]}$, $A_{1}^{2}\subseteq {A_{2}^{1}}^{\left[ \ast \right]}$, $A_{1}^{1}\subseteq {A_{1}^{1}}^{\left[ \ast \right]}$, $A_{2}^{2}\subseteq {A_{2}^{2}}^{\left[ \ast \right]}$.

\item If $A$ is symmetric and $D(A)\cap \mathcal{K}_{i}$ is dense in $\mathcal{K}_{i}$, then $A_{i}^{i}$ is single-valued relation and $A\left( 0 \right)\subseteq \mathcal{K}_{j}$, i.e. $A\left( 0 \right)=A_{i}^{j}\left( 0 \right)=A_{j}^{j}\left( 0 \right), \thinspace j\ne i,\thinspace i,\thinspace j=1,\thinspace 2$.
\end{enumerate}
\end{lemma}

\textbf{Proof:} Note that from $E_{i}(D(A))\subseteq D(A)$ it follows 
$E_{j}\left( D\left( A \right) \right)\subseteq D\left( A \right),\thinspace 
i\ne j$, and
\[
E_{i}(D(A))= \mathcal{K}_{i} \cap D(A), \thinspace i=1, 2.
\]
Then the first two statements of the lemma follow directly from the definition of the relations $A_{i}^{j}$ 

(iii) If $B$ is a linear relation in a Krein space, not necessarily closed, then it holds, 
\[
B\left( 0 \right)\subseteq {D\left( B^{[\ast ]} \right)}^{[\bot 
]}.
\]
Indeed, $y\in B\left( 0 \right)\Rightarrow \left( {\begin{array}{*{20}c}
0\\
y\\
\end{array} } \right)\in B\Rightarrow \left[ y,k \right]=0,\forall \left( 
{\begin{array}{*{20}c}
k\\
k^{'}\\
\end{array} } \right)\in B^{\left[ \ast \right]}\Rightarrow B\left( 0 \right)\subseteq {D\left( B^{[\ast ]} \right)}^{[\bot ]}$.

To prove the converse inclusion $(\supseteq)$ we need assumption that $B$ is closed. Then we have
\[
y\in {D\left( B^{[\ast ]} \right)}^{[\bot ]}\Rightarrow \left[ y,k 
\right]=0,\forall \left( {\begin{array}{*{20}c}
k\\
k^{'}\\
\end{array} } \right)\in B^{\left[ \ast \right]}\Rightarrow \left( 
{\begin{array}{*{20}c}
0\\
y\\
\end{array} } \right)\in {B^{\left[ \ast \right]}}^{\left[ \ast 
\right]}=\bar{B}=B\Rightarrow y\in B\left( 0 \right).
\]
Hence, the converse inclusion holds too, which completes the proof of (\ref{eq22}).

(iv) Let us here clarify notation that we will frequently use in this lemma and the next theorem. For $h_{i}\in D(A)\cap \mathcal{K}_{i}$ it holds
\[
\left( {\begin{array}{*{20}c}
h_{i}\\
h_{i}^{'}\\
\end{array} } \right)\in A\Leftrightarrow \left( {\begin{array}{*{20}c}
h_{i}\\
h_{i}^{'}\\
\end{array} } \right)=\left( {\begin{array}{*{20}c}
h_{i}\\
h_{i}^{1}\left[ + \right]h_{i}^{2}\\
\end{array} } \right)\in A_{i}^{1}+A_{i}^{2},\thinspace i=1,\thinspace 2,
\]
where $h_{i}^{'}=h_{i}^{1}\left[ + \right]h_{i}^{2}\in \mathcal{K}_{1}\left[ + \right]\mathcal{K}_{2}$ and $A_{i}^{1}+A_{i}^{2}$ is operator-like sum. For $h=h_{1}\left[ + \right]h_{2}$ it holds
\[
\left( {\begin{array}{*{20}c}
h\\
h^{'}\\
\end{array} } \right)\in A\Leftrightarrow \thinspace \left( 
{\begin{array}{*{20}c}
h\\
h^{'}\\
\end{array} } \right)=\left( {\begin{array}{*{20}c}
h_{1}\thinspace \thinspace \thinspace \thinspace \thinspace \thinspace 
\left[ + \right]\thinspace \thinspace \thinspace \thinspace \thinspace \thinspace h_{2}\\
h_{1}^{1}\left[ + \right]h_{1}^{2}+h_{2}^{1}\left[ + \right]h_{2}^{2}\\
\end{array} }\thinspace \right), \left( {\begin{array}{*{20}c}
h_{i}\\
h_{i}^{j}\\
\end{array} } \right)\in A_{i}^{j},\thinspace i,\thinspace j=1,\thinspace 
2.
\]
In the sequel we will for addition of vectors frequently use simply $+$ rather than $\left[ + \right]$ because the notation of the vectors in the particular sums indicate when the direct orthogonal sum applies. 

Let us now assume that $A$ is a symmetric relation and let us, for example, show that it holds 
\[
A_{2}^{1}\subseteq {A_{1}^{2}}^{[\ast ]}.
\]
Let us select arbitrary $\left( {\begin{array}{*{20}c}
h_{2}\\
h_{2}^{1}\\
\end{array} } \right)\in A_{2}^{1}$. Then for every 
$\left( {\begin{array}{*{20}c}
h_{1}\\
h_{1}^{2}\\
\end{array} } \right)\in A_{1}^{2}$, 
there exists 
$\left( {\begin{array}{*{20}c} h_{2}\\
h_{2}^{1}+h_{2}^{2}\\
\end{array} } \right)\in A$ and
$\left( {\begin{array}{*{20}c}
h_{1}\\
h_{1}^{1}+h_{1}^{2}\\
\end{array} } \right)\in A$. Because $A$ is symmetric, it holds
\[
\left[ h_{1},h_{2}^{1}+ h_{2}^{2} \right]=\left[ h_{1}^{1}+h_{1}^{2},h_{2} \right].
\]
Hence,
\[
\left[ h_{1},h_{2}^{1} \right]=\left[ h_{1}^{2},h_{2} \right].
\]
This proves $\left( {\begin{array}{*{20}c}
h_{2}\\
h_{2}^{1}\\
\end{array} } \right) \in {A_{1}^{2}}^{[\ast ]}$, i.e. $A_{2}^{1}\subseteq \thinspace {A_{1}^{2}}^{[\ast ]}$. 

By the same token it holds: 
\[
A_{i}^{j}\subseteq {A_{j}^{i}}^{\left[ \ast \right]}, \thinspace \thinspace \forall i, j = 1, 2.
\]

(v) We will prove this statement for $i=1, \thinspace j=2$. Hence, we assume that $D(A)\cap \mathcal{K}_{1}$ is dense in $\mathcal{K}_{1}$. Let us apply formula (\ref{eq22}) on the (closed) relation $B={A_{1}^{1}}^{\left[ \ast \right]}$. We get $B\left( 0 \right)={A_{1}^{1}}^{\left[ \ast \right]}\left( 0 \right)={D\left(\bar{ A_{1}^{1}} \right)}^{\left[ \bot \right]}={D\left( A_{1}^{1} \right)}^{\left[ \bot \right]}=\lbrace 0 \rbrace $. Then it follows 
\[
A_{1}^{1}\left( 0 \right)\subseteq {A_{1}^{1}}^{\left[ \ast \right]}\left( 0 \right)={D\left( A_{1}^{1} \right)}^{\left[ \bot \right]}=\left\{ 0 \right\}
\]
\[
\Rightarrow A\left( 0 \right)=E_{2}A(0)=A_{1}^{2}\left( 0 \right)=A_{2}^{2}\left( 0 \right)\subseteq \mathcal{K}_{2}. \thinspace \thinspace \thinspace \thinspace \square
\] 

In the following theorem, we give necessary and sufficient conditions for a self-adjoint linear relation in a Krein space to be reduced in the sense of Definition \ref{definition4}. The important statement is (vii). Some of the other listed statements are merely the important steps in the proof of the statement (vii). 

\begin{theorem}\label{theorem2} Assume that $A$ is a self-adjoint linear relation in a Krein space $\left( \mathcal{K},\thinspace [.,.] \right)$, $\mathcal{K}_{1}\subset \mathcal{K}$ is a nontrivial non-degenerate subspace, and $\mathcal{K}_{2}$ is the orthogonal complement of $\mathcal{K}_{1}$ in $\mathcal{K}$, i.e.
\[
\mathcal{K}=\mathcal{K}_{1}\left[ + \right]\mathcal{K}_{2}.
\]
If it holds $E_{1}(D(A))\subseteq D(A)$ and $A(\mathcal{K}_{1} \cap D(A)) \subseteq \mathcal{K}_{1}$, then

\begin{enumerate}[(i)]
\item $A=A_{1}^{1}\hat{+}\left( A_{2}^{1}+A_{2}^{2} \right).$
\item $A_{2}^{2}$ is single-valued self-adjoint relation in $\mathcal{K}_{2}.$
\item $A_{2}^{2}$ and $A_{2}^{1}$ are densely defined operators in $\mathcal{K}_{2}.$
\item $A\left( 0 \right)=A_{1}^{1}\left( 0 \right)={D\left( {A_{1}^{1}}^{[\ast ]} \right)}^{[\bot ]}={D(A)}^{[\bot ]}.$
\item ${A_{2}^{1}}^{[\ast ]}$ is single valued.
\item $A_{1}^{1}={A_{1}^{1}}^{[\ast ]}\Leftrightarrow R(A_{2}^{1})\subseteq A_{1}^{1}(0)$
\item $A=A_{1}^{1}[+]A_{2}^{2},$ if and only if $A_{1}^{1}$ is self-adjoint.
\item If $A(D(A)\cap \mathcal{K}_{1})\subseteq \mathcal{K}_{1}$ is dense in $\mathcal{K}_{1}$, then ${A_{1}^{1}}^{[\ast]}=A_{1}^{1}$ is operator as well.
\end{enumerate}
\end{theorem}

\textbf{Proof:} By assumption $h_{1}^{2}\equiv 0, \forall h_{1}\in \mathcal{K}_{1} \cap D(A)$. Then the statement (i) follows from Lemma \ref{lemma6} (ii). 

(ii) Because, $A(\mathcal{K}_{1} \cap D(A))\subseteq \mathcal{K}_{1}$, it holds $A(0)\subseteq \mathcal{K}_{1}$. Hence, $E_{2}A(0)=A_{2}^{2}(0)=\lbrace0\rbrace$, i.e. $A_{2}^{2}$ is single-valued. Let us now prove that $A_{2}^{2}$ is a self-adjoint operator. Assume that 
$\left( {\begin{array}{*{20}c}
k_{2}\\
k_{2}^{2}\\
\end{array} } \right)\in {A_{2}^{2}}^{[\ast ]}$. 
We will first verify that for every 
$\left( {\begin{array}{*{20}c}
h_{1} + h_{2}\\
h_{1}^{1}+h_{2}^{1}+h_{2}^{2}\\
\end{array} } \right)\in A$ 
it holds 
\[
\left[ k_{2},h_{1}^{1}+h_{2}^{1}+h_{2}^{2} \right]=\left[ k_{2}^{2},h_{1} + h_{2} \right].
\]
This equation is obviously equivalent to 
\[
\left[ k_{2},h_{2}^{2} \right]=\left[ k_{2}^{2},h_{2} \right].
\]
which holds according to assumption 
$\left( {\begin{array}{*{20}c}
k_{2}\\
k_{2}^{2}\\
\end{array} } \right)\in {A_{2}^{2}}^{\left[ \ast \right]}$. 
Therefore, 
\\
$\left( 
{\begin{array}{*{20}c}
k_{2}\\
k_{2}^{2}\\
\end{array} } \right)\in A^{\left[ \ast \right]}=A$. Hence, 
$\left( 
{\begin{array}{*{20}c}
k_{2}\\
k_{2}^{2}\\
\end{array} } \right)\in A_{2}^{2}$. This proves (ii). 

(iii) Because, $A_{2}^{2}$ is self-adjoint and single-valued it holds 
\[
\left\{ 0 \right\}=A_{2}^{2}\left( 0 \right)={D\left( {A_{2}^{2}}^{[\ast ]} \right)}^{[\bot ]}.
\]
Hence, $D\left( {A_{2}^{2}}^{\left[ \ast \right]} \right)=D\left( A_{2}^{2} \right)$ is dense in $\mathcal{K}_{2}$. Then also $D\left( A_{2}^{1} \right)=E_{2}\left( D\left( A \right) \right)$ is dense in $\mathcal{K}_{2}$. 

(iv) Because $A$ is self-adjoint and $ A_{1}^{1} \subseteq A $, the following implications hold: 
\[
A_{1}^{1}\subseteq A \Rightarrow A \subseteq {A_{1}^{1}}^{[\ast]} \Rightarrow {D\left( {A_{1}^{1}}^{[\ast]} \right)}^{[\bot ]}\subseteq {D\left( A \right)}^{[\bot ]}.
\]
It also holds $ A_{1}^{1}(0)\subseteq {D\left( {A_{1}^{1}}^{[\ast]} \right)}^{[\bot ]}$, see the proof of (\ref{eq22}). Because, $A=A^{[\ast]}$ is closed we can apply formula (\ref{eq22}) to $A$. We get:
\[
A_{1}^{1}(0)\subseteq {D\left( {A_{1}^{1}}^{[\ast]} \right)}^{[\bot ]}\subseteq {D\left( A \right)}^{[\bot ]}=A(0).
\]
According to the assumption $A(\mathcal{K}_{1} \cap D(A))\subseteq \mathcal{K}_{1}$, it holds $ A(0)=A_{1}^{1}(0)$ and, therefore, the "$\subseteq $" signs become "$=$" signs in the above line, which proves (iv).

(v) ${A_{2}^{1}}^{[\ast ]}\left( 0 \right)={D\left( {A_{2}^{1}}^{\left[ \ast 
\right][\ast ]} \right)}^{[\bot ]}={D\left( A_{2}^{1} \right)}^{[\bot ]}.$ 
According to (iii) $D\left( A_{2}^{1} \right)=D\left( A_{2}^{2} \right)$ is dense in $\mathcal{K}_{2}$. Therefore, 
${A_{2}^{1}}^{[\ast ]}\left( 0 \right)=\left\{ 0 \right\}$, which proves (v).

(vi) Let us first prove
\[
A_{1}^{1}={A_{1}^{1}}^{[\ast ]}\Leftrightarrow R\left( A_{2}^{1} \right)\subseteq A\left( 0 \right).
\]
($\Rightarrow )$ Let us assume that ${A_{1}^{1}}^{[\ast ]}=A_{1}^{1}$, and observe two arbitrary elements
\[
\left( {\begin{array}{*{20}c}
h_{1} + h_{2}\\
h_{1}^{1}+h_{2}^{1}+h_{2}^{2}\\
\end{array} } \right)\in A, \left( 
{\begin{array}{*{20}c}
k_{1} + k_{2}\\
k_{1}^{1}+k_{2}^{1}+k_{2}^{2}\\
\end{array} } \right)\in A.
\]
Because $A$ is self-adjoint, it holds
\[
\left[ h_{1}+h_{2},k_{1}^{1}+k_{2}^{1}+k_{2}^{2} \right]=\left[h_{1}^{1}+h_{2}^{1}+h_{2}^{2},k_{1}+k_{2} \right]\Leftrightarrow 
\]
\[
\left[ h_{1},k_{1}^{1}+k_{2}^{1} \right]+\left[ h_{2},k_{2}^{2} \right]=\left[ h_{1}^{1}+h_{2}^{1},k_{1} 
\right]+\left[ h_{2}^{2},k_{2} \right].
\]
Because,$ A_{1}^{1}$ and $A_{2}^{2}$ are symmetric this equation reduces to 
\[
\left[ h_{1},k_{2}^{1}\thinspace \right]=\left[ h_{2}^{1},k_{1} \right].
\]
Because of $A \left(D \left(A \right) \cap \mathcal{K}_{1} \right)\subseteq \mathcal{K}_{1}$, we have $h_{1}^{2}\equiv 0$.  
Then, according to claim $A_{2}^{1}\subseteq {A_{1}^{2}}^{\left[* \right]}$ in Lemma \ref{lemma6}, it holds 
\[
0=\left[ h_{1}^{2},k_{2} \right]=\left[ h_{1},k_{2}^{1} \right]=\left[ h_{2}^{1},k_{1} \right].
\]
Hence, $R\left( A_{2}^{1} \right)\subseteq {D\left( A_{1}^{1} \right)}^{\left[ \bot \right]}={A_{1}^{1}}^{\left[ \ast \right]}\left( 0 \right)=A_{1}^{1}\left( 0 \right)=A\left( 0 \right)$.

$(\Longleftarrow)$ Assume now that $R\left( A_{2}^{1} \right)\subseteq A\left( 0 \right)$ and prove that $A_{1}^{1}$ is self-adjoint. 

Assume that $\left( {\begin{array}{*{20}c}
k_{1}\\
k_{1}^{1}\\
\end{array} } \right)\in {A_{1}^{1}}^{[\ast ]}$. 
We will first prove that for every 
\[
\left( {\begin{array}{*{20}c}
h_{1} + h_{2}\\
h_{1}^{1}+h_{2}^{1}+h_{2}^{2}\\
\end{array} } \right)\in A
\]
it holds
\begin{equation}
\label{eq24}
\left[ h_{1}+h_{2},k_{1}^{1} \right]=\left[ h_{1}^{1}+h_{2}^{1}+h_{2}^{2},k_{1} \right].
\end{equation}
This equation is equivalent to  
\[
\left[ h_{1},k_{1}^{1} \right]=\left[ h_{1}^{1}+h_{2}^{1},k_{1} \right].
\]
According to our assumption $\left( {\begin{array}{*{20}c}
k_{1}\\
k_{1}^{1}\\
\end{array} } \right)\in {A_{1}^{1}}^{[\ast ]}$, it holds $\left[ 
h_{1},k_{1}^{1} \right]=\left[ h_{1}^{1},k_{1} \right]$. 

It remains to prove $0=\left[ h_{2}^{1},k_{1} \right]$. 
\\
According to our assumption, and (iv), it holds $R\left( A_{2}^{1} \right)\subseteq A\left( 0 
\right)=A_{1}^{1}\left( 0 \right)={D\left( {A_{1}^{1}}^{\left[ \ast \right]} 
\right)}^{\left[ \bot \right]}$. Then we have
\[
R\left( A_{2}^{1} \right)\left[ \bot \right]D\left( {A_{1}^{1}}^{\left[ \ast 
\right]} \right)\Rightarrow \left[ h_{2}^{1},k_{1} 
\right]=0.
\]
Hence, (\ref{eq24}) is satisfied. It further means $\left( 
{\begin{array}{*{20}c}
k_{1}\\
k_{1}^{1}\\
\end{array} } \right)\in A^{\left[ \ast \right]}=A\Rightarrow \left( 
{\begin{array}{*{20}c}
k_{1}\\
k_{1}^{1}\\
\end{array} } \right)\in A_{1}^{1}$. This proves that $A_{1}^{1}$ is 
self-adjoint relation, i. e. it proves $(\Longleftarrow)$. 

Now (vi) follows from $A(0)=A_{1}^{1}(0)$.

(vii)  Assume that $A_{1}^{1}={A_{1}^{1}}^{\left[ \ast \right]}$. According to (i) we have
\[
A=A_{1}^{1}\hat{+}\left( A_{2}^{1}+A_{2}^{2} \right).
\]
According to (vi) we have
\[
A_{1}^{1}={A_{1}^{1}}^{[\ast ]}\Leftrightarrow R(A_{2}^{1})\subseteq A_{1}^{1}(0) \Leftrightarrow
\]
\[
\Leftrightarrow 
\left( {\begin{array}{*{20}c}
 h_{2}\\
h_{2}^{1}+h_{2}^{2}\\
\end{array} } \right)=
\left( {\begin{array}{*{20}c}
 0\\
h_{2}^{1}\\
\end{array} } \right)+
\left( {\begin{array}{*{20}c}
 h_{2}\\
h_{2}^{2}\\
\end{array} } \right), \forall 
\left( {\begin{array}{*{20}c}
 h_{2}\\
h_{2}^{1}+h_{2}^{2}\\
\end{array} } \right) \in A_{2}^{1}+A_{2}^{2}.
\] 
Therefore, for arbitrarily selected element from $A=A_{1}^{1}\hat{+}\left( A_{2}^{1}+A_{2}^{2} \right)$ it holds 
\[
\left( {\begin{array}{*{20}c}
 h_{1}\\
h_{1}^{1}\\
\end{array} } \right)+
\left( {\begin{array}{*{20}c}
 h_{2}\\
h_{2}^{1}+h_{2}^{2}\\
\end{array} } \right)=
\left( {\begin{array}{*{20}c}
 h_{1}\\
h_{1}^{1}\\
\end{array} } \right)+\left( {\begin{array}{*{20}c}
 0\\
h_{2}^{1}\\
\end{array} } \right)+
\left( {\begin{array}{*{20}c}
 h_{2}\\
h_{2}^{2}\\
\end{array} } \right).
\] 
From $A_{2}^{1}\left( 0 \right)\subseteq A_{1}^{1}\left( 0 \right)$ it follows
\[
\left( {\begin{array}{*{20}c}
 h_{1}\\
h_{1}^{1}\\
\end{array} } \right)+\left( {\begin{array}{*{20}c}
 0\\
h_{2}^{1}\\
\end{array} } \right) \in A_{1}^{1}.
\]
Therefore $A=A_{1}^{1} \hat+ A_{2}^{2}$. Because of $[h_{1},h_{2}]=0$ and $[h_{1}^{1}+h_{2}^{1},h_{2}^{2}]=0$ we conclude $A=A_{1}^{1} [+] A_{2}^{2}$.
\\

Conversely, from 
\[
A=A_{1}^{1}[+]A_{2}^{2}
\]
and from Lemma \ref{lemma4} it follows that relations $A_{i}^{i}$ are self-adjoint in the corresponding $\mathcal{K}_{i}, i=1, 2$. 

(viii) This statement also follows from (\ref{eq22}).$\square$

\section{Direct sum representation of generalized Nevanlinna functions}\label{s6}
\textbf{3.1} Let us assume that functions $Q_{i}\in N_{\kappa_{i}}(\mathcal{H})$ are minimally 
represented by triplets $\left( \mathcal{K}_{i},A_{i},\Gamma_{i} \right), i=1, 2$, in representations of the form (\ref{eq2}), where $A_{i}$ are self-adjoint relations in Pontryagin spaces $\mathcal{K}_{i}$, and $\Gamma 
_{i}:\mathcal{H}\to \mathcal{K}_{i}$ are operators. We define the domain of $Q:=Q_{1}+Q_{2}$ by
\[
\mathcal{D}\left( Q \right)=\mathcal{D}\left( Q_{1} \right)\cap \mathcal{D}\left( Q_{2} \right),
\]
space $\tilde{\mathcal{K}}$ as the orthogonal direct sum,
\begin{equation}
\label{eq32}
\tilde{\mathcal{K}}:=\mathcal{K}_{1}\left[ + \right]\mathcal{K}_{2}=c.l.s.\left\{ \left( 
{\begin{array}{*{20}c}
\Gamma_{1z_{1}}h_{1}\\
\Gamma_{2z_{2}}h_{2}\\
\end{array} } \right):z_{i}\in \mathcal{D}\left( Q_{i} \right),h_{i}\in \mathcal{H},\thinspace 
i=1,\thinspace 2\thinspace \right\}.
\end{equation}
Scalar product in $\tilde{\mathcal{K}} $ is naturally defined by
\[
\left[ \left( 
{\begin{array}{*{20}c}
f_{1}\\
f_{2}\\
\end{array} } \right),
\left( 
{\begin{array}{*{20}c}
g_{1}\\
g_{2}\\
\end{array} } \right) \right]
:= \left[ f_{1},g_{1} \right]+\left[ f_{2},g_{2} \right]; \thinspace f_{i}, g_{i} \in \mathcal{K}_{i}, i=1, 2.
\]
In this subsection we will create a minimal state space of $Q$ within $\tilde{\mathcal{K}}$ by means of the elements 
\[
\tilde{\Gamma }_{z}h:=\left( {\begin{array}{*{20}c}
\Gamma_{1z}h\\
\Gamma_{2z}h\\
\end{array} } \right),\thinspace z\in \mathcal{D}\left( Q \right), h\in \mathcal{H}.
\]
First, we will find state manifold of $L \left(Q \right) $. We start with linear space   
\begin{equation}
\label{eq34}
L:=l.s.\left\{ \tilde{\Gamma }_{z}h:z\in \mathcal{D}\left( Q \right),\thinspace h\in \mathcal{H} 
\right\}\subseteq \tilde{\mathcal{K}}.
\end{equation}

The closure of $L$ in $\tilde{\mathcal{K}}$ is given by
\[
\bar{L}={L^{\left[ \bot \right]}}^{\left[ \bot \right]}=c.l.s. \left\{ \tilde{\Gamma }_{z}h:z\in \mathcal{D}\left( Q \right),h\in \mathcal{H} \right\}, 
\]
where $L^{\left[ \bot \right]}$ denotes the orthogonal complement of $L$ in $(\tilde{\mathcal{K}}, \left[ .,. \right])$. It is important to note that, \textit{in general case, an indefinite scalar product} $\left[ .,. \right]$ \textit{may degenerate on the closure of a manifold even if it does not degenerate on the given manifold}, see \cite[p. 39] {IKL}. Later, we will prove that it is not the case with $L$ and $\bar{L}$, see Lemma \ref{lemma10}.

We define operator $\tilde{\Gamma }=\left( {\begin{array}{*{20}c}
\Gamma_{1}\\
\Gamma_{2}\\
\end{array} } \right):\mathcal{H}\to \mathcal{K}_{1}\left[ + \right]\mathcal{K}_{2}$ by 
\[
\tilde{\Gamma }h:=\Gamma_{1}h\left[ + \right]\Gamma_{2}h, \Gamma 
_{i}h\in \mathcal{K}_{i},\thinspace i=1,\thinspace 2.
\]
It holds 
\[
\left[ \tilde{\Gamma }h,k_{1}\left[ + \right]k_{2} \right]=\left( h,\Gamma_{1}^{+}k_{1}+\Gamma_{2}^{+}k_{2} \right),\thinspace \forall k_{1}\left[ + \right]k_{2}\in \tilde{\mathcal{K}}.
\]
Therefore, $\tilde{\Gamma }^{+}:\mathcal{K}_{1}\left[ + \right]\mathcal{K}_{2}\to \mathcal{H}$ satisfies
\[
{\tilde{\Gamma }^{+}=\Gamma }_{1}^{+}+\Gamma_{2}^{+},
\]
where we consider that $\Gamma_{l}^{+}, l=1, 2$, is extended on the whole space 
$\mathcal{K}_{1}\left[ + \right]\mathcal{K}_{2}$ by $\Gamma_{i}^{+}\left( k_{j} \right)=0,\forall k_{j}\in \mathcal{K}_{j},\thinspace \thinspace j\ne i$. 

Let the functions $Q_{i}$ again be minimally represented by (\ref{eq2}). For the function $Q:=Q_{1}+Q_{2}$, consider the following representation
\[
Q\left( z \right)={Q_{1}\left( z_{0} \right)}^{\ast }+{Q_{2}\left( z_{0} \right)}^{\ast }+
\]
\begin{equation}
\label{eq36}
\left( z-\bar{z_{0}} \right)\left( \Gamma_{1}^{+}\thinspace \Gamma_{2}^{+} 
\right)\left( \begin{array}{*{20}c}
I_{1}+\left( z-z_{0} \right){\left( A_{1}-z \right)}^{-1} & 0\\
0 & I_{2}+\left( z-z_{0} \right){\left( A_{2}-z \right)}^{-1}\\
\end{array} \right)\left( {\begin{array}{*{20}c}
\Gamma_{1}\\
\Gamma_{2}\\
\end{array} } \right),
\end{equation}
where $z\in \mathcal{D}\left( Q \right)$ and $I_{i}$ denote identities in $\mathcal{K}_{i}$. Note that (\ref{eq36}) is defined only when $\Gamma_{1}$ and $\Gamma_{2}$ simultaneously map the same vector $h\in \mathcal{H}$ into $\tilde{\mathcal{K}}$. That means that manifold $L$ is the linear span of the vectors 
\begin{equation}
\label{eq38}
\tilde{\Gamma }_{z}h=\left( I+\left( z-z_{0} \right)\left( \tilde{A}-z \right)^{-1} \right)\tilde{\Gamma }h, \thinspace z\in \mathcal{D}\left( Q \right), \thinspace h\in \mathcal{H},
\end{equation}
where the resolvent is defined by 
\[
\left( \tilde{A}-z \right)^{-1}:=\left( {\begin{array}{*{20}c}
\left( A_{1}-z \right)^{-1} & 0\\
0 & \left( A_{2}-z \right)^{-1}\\
\end{array} } \right).
\]
We know that the following holds
\[
\left( \frac{Q_{i}\left( z \right)-Q_{i}\left( \bar{w} \right)}{z-\bar{w}}h_{z},h_{w} \right)=\left[ \Gamma_{iz}h_{z},\Gamma_{iw}h_{w} \right]; z, w\in \mathcal{D}\left( Q_{i} \right), z\ne \bar{w};h_{z}, h_{w} \in \mathcal{H};
\]
\[
\left( {Q_{i}}^{'}\left( z \right)h_{z},h_{\bar{z}} \right)=\left[ \Gamma_{iz}h_{z},\Gamma_{i\bar{z}}h_{\bar{z}} \right], \thinspace i=1,\thinspace 2.
\]
Then it is easy to verify that for function $Q=Q_{1}+Q_{2}$, the following holds
\[
\left( \frac{Q\left( z \right)-Q\left( \bar{w} \right)}{z-\bar{w}}h_{z},h_{w} \right)=\left[ \tilde{\Gamma 
}_{z}h_{z},\tilde{\Gamma }_{w}h_{w} \right]; z, w\in \mathcal{D} \left( Q \right), \thinspace z\ne \bar{w};\thinspace h_{z}, h_{w} \in \mathcal{H};
\]
\[
\left( Q^{'}\left( z \right)h_{z},h_{\bar{z}} \right)=\left[ \tilde{\Gamma}_{z}h_{z},\tilde{\Gamma}_{\bar{z}}h_{\bar{z}} \right].
\]
According to thse equations we can, as in \cite{BL}, identify building blocks of the state manifold $L(Q)$ with the building blocks of $L\subseteq \tilde{\mathcal{K}}$ defined by (\ref{eq34}). In other words, the following holds
\[
\varepsilon_{z}h=\tilde{\Gamma }_{z}h=\left( I+\left( z-z_{0} \right)\left( \tilde{A}-z \right)^{-1} \right)\tilde{\Gamma }h
\]
and $L=L\left( Q \right)$. 
\\

\textbf{3.2} In Section \ref{s4}, we have proved that relation A can be reduced in the sense of Definition \ref{definition4} if it satisfies conditions of Theorem \ref{theorem2}.  In the following theorem we will describe decomposition of $Q$ in terms of the reducing nontrivial subspaces $\mathcal{K}_{i}$ and reducing relations $A_{i}, i=1, 2$, of the representing relation $A$ of $Q$.
\\

\begin{theorem}\label{theorem3} 

(i) Assume
\begin{enumerate}[(a)]
\item A function $Q\in N_{\kappa }\left( \mathcal{H} \right)$ is minimally represented by (\ref{eq2}) and there exist nondegenerate, nontrivial subspaces $\mathcal{K}_{1}$ and $\mathcal{K}_{2}$ that \textbf{reduce} the representing relation $A$, i.e. $A=A_{1}[+]A_{2}$. Then: 
\item $\exists Q_{i}{\in N}_{\kappa_{i}}\left( \mathcal{H} \right), i=1, 2$, minimally represented by the triplets $\left( \mathcal{K}_{i},A_{i},\Gamma_{i}\thinspace \right)$.
\item $Q\left( z \right)=Q_{1}\left( z \right)+Q_{2}\left( z \right), i=1, 2$.
\item The representation (\ref{eq36}) of $Q$ is minimal, i.e. $\mathcal{K}_{1} [+] \mathcal{K}_{2}$ is the minimal state space of $Q$. 
\end{enumerate}
\begin{enumerate}[(ii)]
\item Conversely, if conditions (b), (c) and (d) are satisfied, then the representation (\ref{eq36}) is of the form (\ref{eq2}), and subspaces $\mathcal{K}_{1},\thinspace \mathcal{K}_{2}$ are reducing subspaces of $\tilde{A}:=A_{1}\left[ + \right]A_{2}$, i.e. (a) holds. 
\end{enumerate}
\begin{enumerate}[(iii)]
\item In that case it holds $\kappa_{1}+\kappa_{2}=\kappa $. 
\end{enumerate}
\end{theorem}

\textbf{Proof: }(i) We know that negative index of the minimal state space $\mathcal{K}$ is equal to $\kappa$, the negative index of $Q$.  Let $\mathcal{K}_{1}$ and $\mathcal{K}_{2}$ be nontrivial non-degenerate subspaces that \textbf{reduce} representing relation $A$. Then $\mathcal{K}=\mathcal{K}_{1}\left[ + \right]\mathcal{K}_{2}$ and $A=A_{1}\left[ + \right]A_{2}$. If $\kappa_{i}$, $0\le \kappa_{i}$, denote negative indexes of $\mathcal{K}_{i},i=1,2$, then obviously $\kappa_{1}+\kappa _{2}=\kappa $. 

Because $A$ is a self-adjoint relation, according to Lemma \ref{lemma4},  $A_{i}$ are also self-adjoint relations in $\mathcal{K}_{i}$. Let $E_{i}:\mathcal{K}\rightarrow \mathcal{K}_{i}$ be orthogonal projections, and $\Gamma_{i}:=E_{i}\circ\Gamma, i=1, 2$. Then the following decompositions hold: 
\[
I+\left( z-z_{0} \right)\left( A-zI \right)^{-1}=
\]
\begin{equation}
\label{eq40}
=\left( \begin{array}{*{20}c}
I_{1}+\left( z-z_{0} \right)\left( A_{1}-z \right)^{-1} & 0\\
0 & I_{2}+\left( z-z_{0} \right)\left( A_{2}-z \right)^{-1}\\
\end{array} \right).
\end{equation}
and
\[
Q\left( z \right)=Q_{1}\left( z \right)+Q_{2}\left( z \right),
\]
where
\begin{equation}
\label{eq42}
Q_{i}\left( z \right):={Q_{i}\left( z_{0} \right)}^{\ast }+\left( 
z-\bar{z_{0}} \right)\Gamma_{i}^{+}\left( {I+\left( z-z_{0} \right)\left( 
A_{i}-zI \right)}^{-1} \right)\Gamma_{i}.
\end{equation}
The constant operators ${Q_{i}\left( z_{0} \right)}^{\ast }$ can be arbitrarily selected as long as ${Q_{1}\left( z_{0} \right)}^{\ast}+{Q_{2}\left( z_{0} \right)}^{\ast }={Q\left( z_{0} \right)}^{\ast }$. 
Hence, the minimal representation (\ref{eq2}) of $Q$ can be expressed as the orthogonal sum representation (\ref{eq36}). This proves (c) and (d).

Because $A_{i}$ are self-adjoint linear relations in the Pontryagin spaces $\mathcal{K}_{i}$, functions (\ref{eq42}) are generalized Nevanlinna functions. From (\ref{eq40}) and from the minimality of representation (\ref{eq2}), minimality of representations (\ref{eq42}) follows. 

Indeed, for $y_{1}[+]y_{2}\in \mathcal{K}_{1}[+]\mathcal{K}_{2}$ minimality of (\ref{eq2}) means 
\[
\left[ \left( {\begin{array}{*{20}c}
y_{1}\\
y_{2}\\
\end{array} } \right),\left( {\begin{array}{*{20}c}
I_{1}+{\left( z-z_{0} \right)\left( A_{1}-z \right)}^{-1} & 0\\
0 & I_{2}+{\left( z-z_{0} \right)\left( A_{1}-z \right)}^{-1}\\
\end{array} } \right)\left( {\begin{array}{*{20}c}
\Gamma_{1}h\\
\Gamma_{2}h\\
\end{array} } \right) \right]=0,
\]
\[
\thinspace \forall z\in \rho \left( A \right),\thinspace \forall h\in 
\mathcal{H}\Rightarrow \left( {\begin{array}{*{20}c}
y_{1}\\
y_{2}\\
\end{array} } \right)=0.
\]
If we keep $y_{2}=0$, we can conclude that $Q_{1}$ is minimally represented 
by $\left( \mathcal{K}_{1},A_{1},\Gamma_{1}\thinspace \right)$. By the same token we 
can conclude that $Q_{2}$ is minimally represented by $\left( 
\mathcal{K}_{2},A_{2},\Gamma_{2}\thinspace \right)$. This further means that negative 
indexes of functions $Q_{i}$ are equal to $\kappa_{i}$, the negative indexes of space $\mathcal{K}_{i}$. Hence, $Q_{i}\in N_{\kappa _{i}}(\mathcal{H}),\thinspace i=1,2$. This proves (b). 

From the equation $\kappa_{1}+\kappa_{2}=\kappa $ established for negative 
indexes of $\mathcal{K}_{i}$ and $\mathcal{K}$, now we can conclude that the same equation holds for 
negative indexes of the functions $Q_{i}$ and $Q$. This proves (iii). 

(ii) Assume now that conditions (b), (c) and (d) are satisfied, where $\tilde{A}:=A_{1}\left[ + \right]A_{2}$ is the representing relation of $Q$. Then subspaces $\mathcal{K}_{i}$ and relations $A_{i}$ satisfy conditions of Definition \ref{definition4}, i.e. they are reducing subspaces and reducing relations of the representing relation $\tilde{A}$ in (\ref{eq36}). Because, $A_{i}$ are self-adjoint relations, according to Lemma \ref{lemma4} the relation $\tilde{A}$ is also self-adjoint. According to assumption (d) and Theorem \ref{theorem1}, the triplet $\left(\tilde{\mathcal{K}}, \tilde{A}, \tilde{\Gamma} \right)$ is uniquely determined (up to isomorphism). Hence, representation (\ref{eq36}) is of the form (\ref{eq2}). This proves statement (a), which completes the proof of (ii). $\square$
\\

If the conditions of Theorem \ref{theorem2} are satisfied, then $A_{2}$ is densely defined (single-valued) self-adjoint operator in $\mathcal{K}_{2}$. In that case function $Q_{2}$ has some nice features at infinity, see e.g. \cite[Satz 1.4]{KL2} for scalar functions. If $A_{2}$ is bounded, see \cite[Corollary 1]{B} for operator valued functions.
\\

\textbf{3.3} As discussed in subsection \ref{s6}.1, vectors $\tilde{\Gamma }_{z}h, \thinspace z\in \mathcal{D}\left( Q \right),\thinspace h\in \mathcal{H}$, are building blocks of the state manifold $L(Q)=L$ of $Q:=Q_{1}+Q_{2}$. Let us now consider the structure of $\tilde{\mathcal{K}}$ introduced by (\ref{eq32}). Denote 
\[
L_{0}:=\bar{L}\cap L^{\left[ \bot \right]}.
\]
Recall that the minimal state space $\mathcal{K}$ of $Q$ is defined as completion of the quotient space $\raise0.7ex\hbox{$L$} \!\mathord{\left/ {\vphantom {L L_{0}^{0}}}\right.\kern-\nulldelimiterspace}\!\lower0.7ex\hbox{$L_{0}^{0}$}$, where 
\[
L_{0}^{0}:=L\cap L^{[\bot ]}
\]
(see Section \ref{s2}.1 or \cite{DLS, KL1} for more details). 
\\
Note, $L_{0}^{0}\subseteq L_{0}$ in general case. We will see in Lemma \ref{lemma10} that in our setting it holds  
\[
L_{0}^{0}=\lbrace 0 \rbrace \Rightarrow L_{0} = \lbrace 0 \rbrace. 
\]
For our purpose, we need to decompose $ \tilde{\mathcal{K}}$ by means of $L_{0}$. Obviously, $L_{0}$ is finite-dimensional because it is isotropic subspace of $\bar{L}\subseteq \tilde {\mathcal{K}}$. According to  \cite[Theorem 3.3 and Theorem 3.4]{IKL} the following decompositions hold 
\[
\bar{L}=L_{1}\left[ + \right]L_{0},\thinspace \thinspace L^{[\bot]}=L_{0}\left[ + \right]L_{2},
\]
\begin{equation}
\label{eq44}
\tilde{\mathcal{K}}=L_{1}\left[ + \right]\left( L_{0}\dot{+}F \right)[+]\thinspace L_{2},
\end{equation}
where $L_{1}$ and $L_{2}$ are non-degenerate subspaces and $F$ is a neutral subspace of $\tilde{\mathcal{K}}$, skewly linked to $L_{0}$. Then $\tilde{\kappa}_{0}:=\dim L_{0}$ is the negative index of the non-degenerate subspace $L_{0}\dot{+}F.$ Let $\tilde{\kappa }_{i},\thinspace i=1,\thinspace 2$, denote the negative indexes of subspaces $L_{i}$ in decomposition (\ref{eq44}). $(\mathcal{K},\thinspace A,\thinspace \Gamma )$ again denotes the triplet that minimaly represents $Q=\thinspace Q_{1}+Q_{2}$.

\begin{theorem}\label{theorem6} Let functions $Q_{i}{\in N}_{\kappa_{i}}\left( \mathcal{H} \right)$ be minimally represented by formulas of the form (\ref{eq2}). Assume that the function $Q:=\thinspace Q_{1}+Q_{2}{\in N}_{\kappa }\left( \mathcal{H} \right)$ is represented by orthogonal sum representation (\ref{eq36}).

Then the subspace $L_{1}$ in decomposition (\ref{eq44}) is unitarily equivalent to the minimal state space $\mathcal{K}$ of the function $Q=\thinspace Q_{1}+Q_{2}$. Therefore, $\mathcal{K}$ and $L_{1}$, including the corresponding scalar products, can be identified, i.e. $\mathcal{K}=L_{1}$ and $\tilde{\kappa }_{1}=\kappa $.

\end{theorem}

\textbf{Proof:} Observe representation (\ref{eq36}) of $Q$
\[
Q\left( z \right):={Q\left( z_{0} \right)}^{\ast }+\left( z-\bar{z_{0}} \right)\tilde{\Gamma }^{+}\left( {I+\left( z-z_{0} \right)\left( \tilde{A}-zI \right)}^{-1} \right)\tilde{\Gamma }
\]
and decomposition (\ref{eq44}) of $\tilde{\mathcal{K}}$. In subsection \ref{s6}.1 we have proved that we can consider $\varepsilon_{z}=\tilde{\Gamma }_{z}$, i.e. we can identify manifold $L$ defined by (\ref{eq34}) with the state manifold $L\left( Q \right)$, the starting manifold  in the building of the minimal state space $\mathcal{K}$ of the given function $Q$. Therefore, we can use the usual construction to obtain the minimal Pontryagin state space $\mathcal{K}$ of $Q$ by means of $\tilde{\Gamma }_{z}$ and $L$. Then we will prove that $\mathcal{K}$ is unitarily equivalent to $L_1$. 

Let us first prove that the minimal space $\mathcal{K}$ of $Q=Q_{1}+Q_{2}$, which is equal to completion of $\raise0.7ex\hbox{$L$} \!\mathord{\left/ {\vphantom {L 
L_{0}}}\right.\kern-\nulldelimiterspace}\!\lower0.7ex\hbox{$L_{0}^{0}$}$, is also equal to the completion of  
\[
\raise0.7ex\hbox{$L$} \!\mathord{\left/ {\vphantom {L 
L_{0}}}\right.\kern-\nulldelimiterspace}\!\lower0.7ex\hbox{$L_{0}$}.
\]
For that purpose, let us prove that the naturally defined mapping 
\begin{equation}
\label{eq46}
f+L_{0}\to f+L_{0}^{0},\thinspace \forall f\in L
\end{equation}
is an isometric bijection between $\raise0.7ex\hbox{$L$} \!\mathord{\left/ {\vphantom {L 
L_{0}^{0}}}\right.\kern-\nulldelimiterspace}\!\lower0.7ex\hbox{$L_{0}^{0}$}$ 
and $\raise0.7ex\hbox{$L$} \!\mathord{\left/ {\vphantom {L 
L_{0}}}\right.\kern-\nulldelimiterspace}\!\lower0.7ex\hbox{$L_{0}$}$. 

It obviously holds $L_{0}^{0}\subseteq L_{0}$. Now we have: 
\[
0\ne f+L_{0}\in \raise0.7ex\hbox{$L$} \!\mathord{\left/ {\vphantom {L L_{0}}}\right.\kern-\nulldelimiterspace}\!\lower0.7ex\hbox{$L_{0}$}\Rightarrow f\notin L_{0}^{0}\Rightarrow 0\ne f+L_{0}^{0}\in \raise0.7ex\hbox{$L$} \!\mathord{\left/ {\vphantom {L L_{0}^{0}}}\right.\kern-\nulldelimiterspace}\!\lower0.7ex\hbox{$L_{0}^{0}$}.
\]
In order to prove the converse implication, let us assume the contrary:
\[
0\ne f+L_{0}^{0}\in \raise0.7ex\hbox{$L$} \!\mathord{\left/ {\vphantom {L 
L_{0}^{0}}}\right.\kern-\nulldelimiterspace}\!\lower0.7ex\hbox{$L_{0}^{0}$}\thinspace \thinspace \thinspace and \thinspace \thinspace \thinspace 0=f+L_{0}\in \raise0.7ex\hbox{$L$} \!\mathord{\left/ {\vphantom 
{L L_{0}}}\right.\kern-\nulldelimiterspace}\!\lower0.7ex\hbox{$L_{0}$}.
\]
Then $0=f+L_{0}$ means $f\in L\thinspace $ and $\thinspace f\in L_{0}$. It further means 
that $f\in L$ and $f[\bot ]\bar{L}$. Because, $\bar{L}\supseteq L$ it 
follows $f[\bot ]L$; hence $f\in L_{0}^{0}$, which is a contradiction. This 
proves that the naturally defined mapping (\ref{eq46}) is a bijection, and we can 
identify $\raise0.7ex\hbox{$L$} \!\mathord{\left/ {\vphantom {L 
L_{0}^{0}}}\right.\kern-\nulldelimiterspace}\!\lower0.7ex\hbox{$L_{0}^{0}$}$ 
and $\raise0.7ex\hbox{$L$} \!\mathord{\left/ {\vphantom {L 
L_{0}}}\right.\kern-\nulldelimiterspace}\!\lower0.7ex\hbox{$L_{0}$}$.

Recall that the scalar product is introduced in $\mathcal{K}$ in the following manner: 
If for $f,\thinspace g,\thinspace ...\in L$ the corresponding classes in the 
quotient space $\raise0.7ex\hbox{$L$} 
\!\mathord{\left/ {\vphantom {L 
L_{0}^{0}}}\right.\kern-\nulldelimiterspace}\!\lower0.7ex\hbox{$L_{0}^{0}$}$ are denoted by $\hat{f},\thinspace \hat{g},\thinspace 
\mathellipsis $, then the scalar product is defined by 
\begin{equation}
\label{eq48}
\langle\hat{f},\hat{g}\rangle:=\left[f,g \right].
\end{equation}
Then the quotient space $\raise0.7ex\hbox{$L$} \!\mathord{\left/ {\vphantom 
{L 
L_{0}}}\right.\kern-\nulldelimiterspace}\!\lower0.7ex\hbox{$L_{0}$}=\raise0.7ex\hbox{$L$} 
\!\mathord{\left/ {\vphantom {L 
L_{0}^{0}}}\right.\kern-\nulldelimiterspace}\!\lower0.7ex\hbox{$L_{0}^{0}$}$ 
can be completed in the usual way, see e.g. [5, Section 2.4]. The completion 
$\mathcal{K}$ is unitarily equivalent to space $L_{1}$ introduced by equations (\ref{eq44}). 

Indeed, according to the above definition (\ref{eq48}) and \cite[Theorem 2.4 (i)]{IKL}, the 
sequence $\left\{ f_{n} \right\}\subseteq L_{1}\cap L$ converges to some 
$f_{0}\in L_{1}$ if and only if the sequence $\left\{ \hat{f_{n}} 
\right\}_{n=1}^{\infty }=\left\{ f_{n}+L_{0} \right\}_{n=1}^{\infty 
}\subseteq \mathcal{K}$ converges to $\hat{f_{0}}\in \mathcal{K}$. Therefore, equation (\ref{eq48}) 
extends to $L_{1}$ and $\mathcal{K}$. This proves that $\mathcal{K}$ and $L_{1}$ are unitarily 
equivalent and we can consider 
\[
\mathcal{K}=L_{1}
\]
and $\tilde{\kappa }_{1}=\kappa $. $\square$

\begin{remark}\label{remark5} From (\ref{eq38}) it follows 
\[
\tilde{\Gamma }_{z_{0}}=\tilde{\Gamma }
\]
and, therefore 
\[
L_{0}\subset L^{[\bot ]}\subseteq \ker \tilde{\Gamma }^{+}.
\]
Hence, the operator $\Gamma^{+}:\mathcal{K}\to \mathcal{H}$ defined by $\Gamma 
^{+}\hat{f}:=\tilde{\Gamma }^{+}f$ is well defined. If we also set $\Gamma 
h:=\tilde{\Gamma }h,\thinspace \forall h\in \mathcal{H}$, $A:=\tilde{A}_{\mid L_{1}}$ in 
representation (\ref{eq36}) of $Q$, then we obtain representation (\ref{eq2}). 
\end{remark}

\begin{corollary}\label{corollary5}Let functions $Q_{i}{\in N}_{\kappa_{i}}\left( \mathcal{H} \right)$ be minimally represented by formulas of the form (\ref{eq2}) and $Q:=\thinspace Q_{1}+Q_{2}$.The following statements hold:

\begin{enumerate}[(i)]
\item $\tilde{\mathcal{K}}$ is the minimal state space of $Q$ if and only if $L_{1}=\bar{L}=\tilde{\mathcal{K}}$. In that case $\kappa =\kappa_{1}+\kappa_{2}$.
\item $\kappa =\kappa_{1}+\kappa_{2}$ if and only if $\tilde{\mathcal{K}}=L_{1}\left[ + \right]\thinspace L_{2}$, where $L_{2}=\lbrace 0\rbrace $ or $L_{2}=L^{[\bot ]}$ is a positive subspace. 
\item $L_{0}=\left\{ 0 \right\}$ is necessary but not sufficient condition for $\kappa =\kappa_{1}+\kappa_{2}$. 
\end{enumerate}
\end{corollary}
\textbf{Proof:} (i) Assume, $\tilde{\mathcal{K}}$ is minimal state space of $Q$. According to first equation of (\ref{eq44}) it holds $L_{1}\subseteq\bar{L}\subseteq\tilde{\mathcal{K}}$. According to Theorem \ref{theorem6}, $L_{1}$ is minimal state space of $Q$.  Therefore, $L_{1}=\bar{L}=\tilde{\mathcal{K}}$. 

Conversely, if $L_{1}=\bar{L}=\tilde{\mathcal{K}}$ holds, then minimality of $\tilde{\mathcal{K}}$ follows from Theorem \ref{theorem6}. Then $\kappa =\kappa_{1}+\kappa_{2}$ follows from Theorem \ref{theorem3}.

(ii)  Assume $\kappa =\kappa_{1}+\kappa_{2}$. That means that the numbers of negative squares of  $L_{1}$ and $\tilde{\mathcal{K}}$ are equal, and $\tilde{\kappa}=0$. According to (\ref{eq44}) it must be $L_{0}=\lbrace 0\rbrace$. Therefore, $\tilde{\mathcal{K}}=L_{1} [+]L_{2}$, where $L_{2}=\lbrace 0\rbrace $ or $L_{2}$ is a positive subspace. 

Conversely, $\tilde{\mathcal{K}}=L_{1} [+]L_{2}$ and $L_{2}=\lbrace 0\rbrace$ or $L_{2}$  is positive, means that the numbers of negative squares of $L_{1}$ and $\tilde{\mathcal{K}}$ are equal, i.e. $\kappa =\kappa_{1}+\kappa_{2}$.  

In Example \ref{example4} we will prove that there exists the case where $\tilde{\mathcal{K}}=L_{1} [+]L_{2}$ and $L_{2} $ is positive.  

(iii)  $\kappa =\kappa_{1}+\kappa_{2} \Rightarrow L_{0} = \lbrace 0 \rbrace \Rightarrow L_{0}^{0}= \lbrace 0 \rbrace $. In Example \ref{example6} we will show that there exists the case where $\tilde{\mathcal{K}}=L_{1} [+]L_{2}$ and $L_{2} $ is negative subspace. That is an example where it holds $L_{0}=\lbrace 0\rbrace$ and $\kappa < \kappa_{1}+\kappa_{2}$. $\square $ 

\section{Analytic criteria}\label{s8}
\textbf{4.1} In this section, we will prove criteria that enable us to research the underlying state space, and negative index of the sum $Q:=Q_{1}+Q_{2}$ analytically, without knowing operator representations of $Q$, $Q_{1}$, $Q_{2}$. In order to derive equations in those criteria we will have to use Definition \ref{definition2} and definitions of scalar products in terms of formal sums, see Section 1.1.  
\\

Let us consider any function $Q \in N_{\kappa}(\mathcal{H})$. By definition $\kappa $ is the maximal (finite) number of negative squares of the sesquilinear form $[.,.]$ defined by the sums 
\begin{equation}
\label{eq48a}
\sum\limits_{i,j=1}^{n} \left[ \Gamma _{z_{i}}h_{i},\Gamma_{z_{j}}h_{j} \right] := \sum\limits_{i,j=1}^{n} \left( \frac{Q\left( z_{i} \right)-Q\left( \bar{z_{j}} \right)}{z_{i}-\bar{z_{j}}}h_{i},h_{j} \right), 
\end{equation}
where $z_{l}\in \mathcal{D}(Q)$, $h_{l}\in \mathcal{H},\thinspace l=1, ..., n$. In other words, $\kappa $ is the negative index of the state manifold $\left( L(Q), [.,.] \right)$. According to Theorem \ref{theorem1}, the \textit{negative index of the \textbf{minimal} state space $\mathcal{K}$ is also equal to $\kappa $}. 

Let us now focus on the sum $Q=Q_{1}+Q_{2}$. Then sum (\ref{eq48a}) can be written as 
\[
\sum\limits_{i,j=1}^{n} \left( \left( \frac{Q_{1}\left( z_{i} \right)-Q_{1}\left( \bar{z_{j}} \right)}{z_{i}-\bar{z_{j}}}+\frac{Q_{2}\left( z_{i} \right)-Q_{2}\left( \bar{z_{j}} \right)}{z_{i}-\bar{z_{j}}} \right)h_{i},h_{j} \right) 
=\sum\limits_{i,j=1}^{n} \left[ \left( 
{\begin{array}{*{20}c}
\Gamma_{1z_{i}}h_{i}\\
\Gamma_{2z_{i}}h_{i}\\
\end{array} } \right),\left( {\begin{array}{*{20}c}
\Gamma_{1z_{j}}h_{j}\\
\Gamma_{2z_{j}}h_{j}\\
\end{array} } \right) \right], 
\] 
where $z_{l}\in \mathcal{D}\left( Q_{1} \right)\cap \mathcal{D}\left( Q_{2} \right)=:\mathcal{D}\left( Q \right), \thinspace h_{l}\in \mathcal{H}, \thinspace l=1, ..., n $. Such sums are subset of sums (\ref{eq49}) below, which generate the inner product in $\tilde{\mathcal{K}} := \mathcal{K}_{1}[+]\mathcal{K}_{2} $. Indeed, here $Q_{1}$ and $Q_{2}$ take the same domain points $z_{l} \in \mathcal{D}(Q)$, while in (\ref{eq49}) $Q_{1}$ and $Q_{2}$ take domain points $z_{l} \in \mathcal{D}(Q_{1})$ and $\zeta_{l} \in \mathcal{D}(Q_{2})$  independently. This means that the space \textit{$\tilde{\mathcal{K}}$ created by means of the sums (\ref{eq49}) may be larger than the state space $\mathcal{K}$, which is created by means of the sums (\ref{eq48a})}.
\\ 

Now we can prove the following lemma. 

\begin{lemma}\label{lemma10} Assume that functions $Q_{i}{\in N}_{\kappa_{i}}\left( \mathcal{H} \right)$ are minimally represented by triplets $\left( \mathcal{K}_{i},A_{i},\Gamma_{i}\thinspace \right),\thinspace i=1,\thinspace 2$, and $Q:=Q_{1}+Q_{2}$. If scalar product does not degenerate on the state manifold $L=L(Q)$, i.e. if $L_{0}^{0}=\lbrace 0 \rbrace$, then scalar product does not degenerate on $\bar{L}$, and it holds $\tilde{\mathcal{K}}=\mathcal{K}[+]L_{2}$ , where $\mathcal{K}=\bar{L}$ is the minimal state space of $Q$.
\end{lemma}

\textbf{Proof:} According to (\ref{eq34}), $L \subseteq \tilde{\mathcal{K}}$. Let us assume that form $[.,.]$ induced by (\ref{eq48a}) in the state manifold $L=L(Q)$ does not degenerate, i.e. $L_{0}^{0}=\lbrace 0 \rbrace$. Then $\raise0.7ex\hbox{$L$} \!\mathord{\left/ {\vphantom {L L_{0}^{0}}}\right.\kern-\nulldelimiterspace}\!\lower0.7ex\hbox{$L_{0}^{0}$}=L $, and \textit{the minimal state space $\mathcal{K}$ is by definition equal to the \textbf{completion} of $L$}.  

Because Pontryagin space $\tilde{\mathcal{K}}$ is complete, the closure $\bar{L}\subseteq \tilde{\mathcal{K}}$ is also complete. Then it holds 
\[
L \subseteq \mathcal{K} \subseteq \bar{L}.
\]
Pontryagin space, $\mathcal{K}$ is non-degenerate.  Because, completion $\mathcal{K}$ is a closed set in $\tilde{\mathcal{K}}$, and $\bar{L}$ is the smallest closed set which contains $L$, we conclude $\mathcal{K}=\bar{L}$. Hence, $\bar{L}$ is non-degenerate. 

Then according to (\ref{eq44}) it holds $\tilde{\mathcal{K}}=\mathcal{K}[+]L_{2}$. $\square$
\\

\textbf{4.2} By definition of $\tilde{\mathcal{K}}$, see (\ref{eq32}), the negative index $\kappa := \kappa _{1}+\kappa_{2}$ of $\tilde{\mathcal{K}}$ is equal to the maximal number of negative squares of the form defined by means of the sums
\[
\sum\limits_{i,j=1}^{n} \left[ \left( {\begin{array}{*{20}c}
\Gamma_{1z_{i}}h_{i}\\
\Gamma_{2\mathbf{\zeta }_{i}}f_{i}\\
\end{array} } \right),\left( {\begin{array}{*{20}c}
\Gamma_{1z_{j}}h_{j}\\
\Gamma_{2\mathbf{\zeta }_{j}}f_{j}\\
\end{array} } \right) \right] =
\]
\begin{equation}
\label{eq49}
\sum\limits_{i,j=1}^{n} {\left( \frac{Q_{1}\left( z_{i} \right)-Q_{1}\left( \bar{z_{j}} \right)}{z_{i}-\bar{z_{j}}}h_{i},h_{j} \right)+\left( \frac{Q_{2}\left( \zeta_{i} \right)-Q_{2}\left( \bar{\zeta 
_{j}} \right)}{\zeta_{i}-\bar{\zeta_{j}}}f_{i},f_{j} \right)},
\end{equation}
where $z_{l}\in \mathcal{D}\left( Q_{1} \right),\thinspace \mathbf{\zeta }_{l}\in \mathcal{D}\left( Q_{2} \right);\thinspace h_{l},\thinspace f_{l}\in \mathcal{H}, \thinspace l=1, ..., n$. Because points $z_{l}, \zeta_{l}$ are arbitrarily selected in their domains, we can create the following sums out of (\ref{eq49}).
\begin{equation}
\label{eq50}
\sum\limits_{i,j=1}^{n} \left[ \left( {\begin{array}{*{20}c}
\Gamma_{1w_{i}}h_{i}\\
\Gamma_{2w_{i}}h_{i}\\
\end{array} } \right),\left( {\begin{array}{*{20}c}
\Gamma_{1w_{j}}h_{j}\\
\Gamma_{2w_{j}}h_{j}\\
\end{array} } \right) \right] 
\mathbf{+}\sum\limits_{i,j=1}^{n} \left[ \left( 
{\begin{array}{*{20}c}
\Gamma_{1z_{i}}h_{i}\\
\Gamma_{2\mathbf{\zeta }_{i}}f_{i}\\
\end{array} } \right),\left( {\begin{array}{*{20}c}
\Gamma_{1z_{j}}h_{j}\\
\Gamma_{2\mathbf{\zeta }_{j}}f_{j}\\
\end{array} } \right) \right], 
\end{equation}
where $w_{l} \in \mathcal{D}\left( Q \right), z_{l} \in \mathcal{D}\left( Q_{1} \right), \zeta_{l} \in \mathcal{D}\left( Q_{2} \right)$, and the second sum is created by vectors that satisfy condition 
\[
\left( {\begin{array}{*{20}c}
\Gamma_{1z_{l}}h_{l}\\
\Gamma_{2\mathbf{\zeta }_{l}}f_{l}\\
\end{array} } \right)\left[ \bot \right]L.
\]
Note that the first sum here is associated with $L$. The  orthogonality condition for vectors from the second sum in (\ref{eq50}) can be written with simplified notation as:
\[
\left[ \left( {\begin{array}{*{20}c}
\Gamma_{1z}h_{1}\\
\Gamma_{2 \zeta}h_{2}\\
\end{array} } \right),\left( {\begin{array}{*{20}c}
\Gamma_{1w}g\\
\Gamma_{2w}g\\
\end{array} } \right) \right]=0,\thinspace \forall w\in \mathcal{D}\left( Q \right),\thinspace \forall g\in \mathcal{H},
\]
where $z\in \mathcal{D}\left( Q_{1} \right),\thinspace \zeta \in \mathcal{D}\left( Q_{2} \right), \thinspace h_{i}\in \mathcal{H}$, $i=1, 2$. Because, scalar product $(.,.)$ in $\mathcal{H}$ is non-degenerate, this condition can be written as the equation 
\begin{equation}
\label{eq52}
\frac{Q_{1}\left( z \right)-Q_{1}\left( \bar{w} 
\right)}{z-\bar{w}}h_{1}+\thinspace \frac{Q_{2}\left( \zeta 
\right)-Q_{2}\left( \bar{w} \right)}{\zeta -\bar{w}}h_{2}=0,\thinspace 
\forall w\in \mathcal{D}\left( Q \right).
\end{equation}

\begin{lemma}\label{lemma12} Let $Q \in N_{\kappa}\left( \mathcal{H} \right)$ be any minimally represented function by a triplet $\left( \mathcal{K},A,\Gamma \right)$. 
\begin{enumerate}[(i)]
\item If there exist $z \in \mathcal{D}(Q)$ such that $\ker{\Gamma_{z}} \neq \lbrace 0 \rbrace$, then 
\[
\ker{\Gamma_{z}}=\ker{\Gamma_{w}}=:\ker{\Gamma}; \thinspace \forall w \in \mathcal{D}(Q).
\]
\item $h\in \ker{\Gamma}$ if and only if 
\[
\frac{Q\left( z \right)-Q\left( \bar{w} 
\right)}{z-\bar{w}}h = 0, \thinspace \forall z, \forall w \in \mathcal{D}(Q).
\]
\end{enumerate}
\end{lemma}

\textbf{Proof:} (i) For function $Q$ minimally represented by (\ref{eq2}), it holds $\rho(A)= \mathcal{D}(Q)$ and 
\[
\Gamma_{z}=\left( I+\left( z-w \right)\left( A-z \right)^{-1} \right)\Gamma_{w}, \forall w \in \rho(A)= \mathcal{D}(Q), 
\]
see \cite{DLS,HSW}. Assume the contrary to the claim (i), that for some $w \in \mathcal{D}(A)$ it holds $\Gamma_{w} h \neq 0, \Gamma_{z} h = 0 $. Then we have
\[
0=\Gamma_{z}h=\left(I+\left(z-w \right)\left(A-z \right)^{-1} \right)\Gamma_{w}h\Rightarrow \left(z-w \right)\left(A-z \right)^{-1}\Gamma_{w}h=-\Gamma_{w}h.
\]
According to \cite[2.11]{A} it holds
\[
\left(A-z \right)\left(A-z \right)^{-1}\supseteq I \Rightarrow  \left(z-w \right)\Gamma_{w}h \subseteq -\left(A-z \right)\Gamma_{w}h.
\]
Therefore $w \Gamma_{w}h \in A\left(\Gamma_{w}h \right)$, i.e.  $w$ is an eigenvalues of $A$. This contradicts to the fact that $w$ is a regular point of $A$. This proves $\ker{\Gamma_{z}}\subseteq\ker{\Gamma_{w}}$. The converse inclusion is obvious. This proves the first equation of (i). 

Because, $\ker{\Gamma_{w}}$ is independent of $w\in \mathcal{D}(Q) $, we can introduce $\ker{\Gamma}:=\ker{\Gamma_{w}}, \\
w\in \mathcal{D}(Q)$. It is obvious now that claim (i) holds for any two points $z, \thinspace w \in \mathcal{D}(Q)$. This completes the proof of (i). 
\\

(ii) We have
\[
h\in \ker{\Gamma} \Leftrightarrow \left[ \Gamma _{z}h,\Gamma_{w}g \right]=0, \forall z\in \mathcal{D}(Q),  \forall w\in \mathcal{D}(Q),\forall g \in \mathcal{H} 
\]
\[
\Leftrightarrow \Gamma_{w} ^{+}\Gamma_{z}h=0, \forall z\in \mathcal{D}(Q),  \forall w\in \mathcal{D}(Q) \Leftrightarrow \frac{Q\left( z \right)-Q\left( \bar{w} \right)}{z-\bar{w}}h = 0, \thinspace \forall z, \forall w \in \mathcal{D}(Q). 
\] 
\\
The following statement is a criteria that identifies zero-symbols $\varepsilon_{z}h=\Gamma_{z}h$, i.e. the symbols that do not play any role in the state manifold $L(Q)$.  

\begin{corollary}\label{corollary6} Let $Q \in N_{\kappa}\left( \mathcal{H} \right)$ be a minimally represented function by a triplet $\left( \mathcal{K},A,\Gamma \right)$. If there exist a solution $(z_{0},h) \in \mathcal{D}(Q)\times \mathcal{H}, \thinspace h \neq 0$ of the equation
\begin{equation}
\label{eq53a}
\frac{Q\left( z \right)-Q\left( \bar{w} \right)}{z-\bar{w}}h = 0, \thinspace \forall w \in \mathcal{D}(Q),
\end{equation}
then $\Gamma_{z}h=0, \forall z \in \mathcal{D}(Q)$.  
\end{corollary}

It is easy to find regular matrix functions that satisfy (\ref{eq53a}), i.e. that have $\ker \Gamma \neq \lbrace 0 \rbrace$.

\begin{example}\label{example2} Consider the following regular matrix functions:
\[
Q\left( z \right)=\left( {\begin{array}{*{20}c}
z+a & z\\
z & z+b\\
\end{array} } \right)\in N_{\kappa} \left( \mathbf{C}^{2} \right);\thinspace a, \thinspace b \in R, \thinspace(a,b)\neq (0,0), \thinspace \kappa \in \lbrace 0, \thinspace 1, \thinspace 2\rbrace .
\]
Then for vector $h=\left( {\begin{array}{*{20}c}
1\\
-1\\
\end{array} } \right)$, identity (\ref{eq53a}) holds.  \thinspace \thinspace \thinspace \thinspace $\square$
\end{example}

Now we can classify solutions of equation (\ref{eq52}). According to Lemma \ref{lemma12}, if \\
$\left( {\begin{array}{*{20}c}
h_{1}\\
h_{2}\\
\end{array} } \right)\in \ker \Gamma_{1} \times \ker \Gamma_{2}$, then for \textbf{both} functions $Q_{i}$ it holds 
\[
\frac{Q_{i}\left( z_{i} \right)-Q_{i}\left( \bar{w} \right)}{z_{i}-\bar{w}}h_{i} \equiv 0.
\]

Let us call such solutions $\left( {\begin{array}{*{20}c}
h_{1}\\
h_{2}\\
\end{array} } \right)$ of (\ref{eq52}) \textit{singular} solutions. Then, according to Lemma \ref{lemma12} the vectors $\left( {\begin{array}{*{20}c}
\Gamma_{1z_{1}}h_{1}\\
\Gamma_{2z_{2}}h_{2}\\
\end{array} } \right)\equiv 0, \thinspace \forall \left( {\begin{array}{*{20}c}
z_{1}\\
z_{2}\\
\end{array} } \right) \in \mathcal{D}(Q_{1})\times \mathcal{D}(Q_{2}), \thinspace i=1, 2$, i.e. \textbf{they do not exist in $\tilde{\mathcal{K}}$}. Therefore, we can exclude singular solutions of (\ref{eq52}) from the following considerations about structure of $\tilde{\mathcal{K}}$, without loss of generality. Hence, in the following definitions we assume that we deal only with \textbf{non-singular} solutions. It is consistent with the standard assumption that the functions $\Gamma_{z}$ are injections. 
\\

The obvious solutions $\left( z_{1},z_{2};h_{1},h_{2} \right)\in \mathcal{D}(Q_{1}) \times \mathcal{D}(Q_{2}) \times \mathcal{H} \times \mathcal{H}$ of (\ref{eq52}), i.e. the solutions with $h_{1}=h_{2}=0$, we call \textit{trivial solutions}. Hence, the \textbf{non-singular} solutions of (\ref{eq52}) with $\left( h_{1},h_{2} \right)\ne (0,0)$, we call \textit{\textbf{nontrivial}}. We will solve equation (\ref{eq52}), later in couple of examples. 

Let us introduce expression
\[
E=E\left( z_{1},z_{2};h_{1},h_{2} \right):=\left( \frac{Q_{1}\left( z_{1} 
\right)-Q_{1}\left( \bar{z_{1}} \right)}{z_{1}-\bar{z_{1}}}h_{1},h_{1} 
\right)+\left( \frac{Q_{2}\left( z_{2} \right)-Q_{2}\left( \bar{z_{2}} 
\right)}{z_{2}-\bar{z_{2}}}h_{2},h_{2} \right).
\]
i.e.
\[
E=E\left( z_{1},z_{2};h_{1},h_{2} \right):=\left[ \left( {\begin{array}{*{20}c}
\Gamma_{1z_{1}}h_{1}\\
\Gamma_{2z_{2}}h_{2}\\
\end{array} } \right),\left( {\begin{array}{*{20}c}
\Gamma_{1z_{1}}h_{1}\\
\Gamma_{2z_{2}}h_{2}\\
\end{array} } \right) \right].
\]
A nontrivial solution $\left( z_{1},z_{2};h_{1},h_{2} \right)$ of (\ref{eq52}) we call\textit{ positive, negative, neutral }if it satisfies 
\[
E\left( z_{1},z_{2};h_{1},h_{2} \right)>0,\thinspace <0,\thinspace =0,\thinspace respectively.
\]
For $z_{1}=z_{2}=z$ and $h_{1}=h_{2}=h$ we get an important special case of equation (\ref{eq52}):
\begin{equation}
\label{eq53}
\left( \frac{Q_{1}\left( z \right)-Q_{1}\left( \bar{w} 
\right)}{z-\bar{w}}+\thinspace \frac{Q_{2}\left( z \right)-Q_{2}\left( 
\bar{w} \right)}{z-\bar{w}} \right)h=0,\thinspace \forall w\in \mathcal{D}\left( Q 
\right).
\end{equation}

\textbf{Why is this equation important?} Equation (\ref{eq53}) identifies when term $\frac{Q_{1}\left( z \right)-Q_{1}\left( \bar{w} \right)}{z-\bar{w}}$ cancels out with term $\frac{Q_{2}\left( z \right)-Q_{2}\left( \bar{w} \right)}{z-\bar{w}}$. That is how a negative square is lost, i.e. the negative index is  reduced in sum $Q_{1}+Q_{2}$. Then in the underlying space $\tilde{\mathcal{K}}$ we have: 

Assume that $(z;h)$ is a nontrivial (and non-singular) solution of (\ref{eq53}). That means that there exists a nonzero vector $ \tilde{\Gamma}_{z}h := \left( {\begin{array}{*{20}c}
\Gamma_{1z}h\\
\Gamma_{2z}h\\
\end{array} } \right) \in \tilde{K}$. On the other hand, according to Corollary \ref{corollary6}, for the symbol $\Gamma_{z}h$ corresponding to $Q(z)$ in the minimal state space $\mathcal{K}$ of $Q$ it holds $\Gamma_{z}h=0$. Hence, we learn that $(0 \neq ) \tilde{\Gamma}_{z}h \in L_{0}^{0} \subset \tilde{\mathcal{K}}$ corresponds to $(0=)\Gamma_{z}h:=\tilde{\Gamma}_{z}h +L_{0}^{0} \in \raise0.7ex\hbox{$L$} \!\mathord{\left/ {\vphantom {L L_{0}^{0}}}\right.\kern-\nulldelimiterspace}\!\lower0.7ex\hbox{$L_{0}^{0}$} \subseteq \mathcal{K}$.

Let us interpret this explanation in terms of almost Pontryagin spaces. Recall that an \textit{almost Pontryagin space } is a Pontryagin space to which a finite dimensional degenerate linear space has been added orthogonaly, see \cite{SW}. 

According to Theorem \ref{theorem6} we have $\bar{L}=L_{1}[+]L_{0}=\mathcal{K}[+]L_{0}$, i.e. $\bar{L}$ is an almost Pontryagin space, with isotropic subspace $L_{0}$. By similar method we obtain almost Pontryagin spaces $\bar{L^{i}}=K_{i}[+]L^{i}_{0}, \thinspace i=1,2$. Because $L^{i}_{0}\cap \tilde{K}=\lbrace 0 \rbrace, \thinspace i=1,2$, the overlap $\bar{L^{1}}\cap \bar{L^{2}}$ does not have any nonzero elements in $ \tilde{\mathcal{K}} $. The symbols $ \Gamma_{z}, \thinspace z\in \mathcal{D}(Q), h \in \mathcal{H} $ that belong to the overlap are characterized as singular solutions of the equation (\ref{eq52}), and excluded from the considerations. Hence, the overlap does not affect the negative index $\kappa$. However, the negative index $\kappa$ is affected by the existence of nonzero elements $\tilde{\Gamma}_{z}h $ in the isotropic subspace $L_{0}$ of the almost Pontryagin space $\bar{L}$. Those elements are characterized by the non-singular, nontrivial solutions of the equation (\ref{eq53}).
\\

The following theorem gives us further analytic means to investigate structure of the state space $\tilde{\mathcal{K}}$ and to compare number of negative squares $\kappa$ of $Q:=Q_{1}+Q_{2}$ with the sum $\kappa_{1}+\kappa_{2}$.

\begin{theorem}\label{theorem8} Assume $Q_{i}\in N_{\kappa_{i}}\left( \mathcal{H} \right)$ are functions minimally represented by triplets $\left( \mathcal{K}_{i},A_{i},\Gamma_{i} \right), i=1, 2,$ and $Q:=Q_{1}+Q_{2}$ is represented by $( \tilde{\mathcal{K}},\tilde{A,}\thinspace \tilde{\Gamma } )$, i.e. by (\ref{eq36}). 

\begin{enumerate}[(i)]
\item There exists a non-trivial solution of equation (\ref{eq52}) if and only if $\tilde{\mathcal{K}}$ is \textbf{not} minimal state space of $Q$.

\item Equation (\ref{eq53}) has a nontrivial solution $(z,h) \in \mathcal{D}(Q)\times \mathcal{H} $, if and only if $\varepsilon_{z}h:=\tilde{\Gamma }_{z}h\in L_{0}^{0}$.

\item If any nontrivial solution of equation (\ref{eq52}) is neutral or negative, then $\kappa <
\kappa_{1}+\kappa_{2}$.

\item If all nontrivial solutions of (\ref{eq52}) are positive, then $L^{\left[ \bot \right]}=L_{0}[+] L_{2}$ is a non-negative subspace of $\tilde{\mathcal{K}}$, where $L_{2}$ is positive definite subspace.  In this case $\kappa = \kappa_{1}+\kappa_{2}$ if and only if the state manifold $L=L(Q)$ is non-degenerate.

\item A necessary condition for $\kappa =\thinspace \kappa_{1}+\kappa_{2}$ is that equation (\ref{eq53}) has only trivial solutions. Evan the stronger condition, that equation (\ref{eq52}) has only trivial solutions, is not sufficient for $\kappa =\thinspace \kappa_{1}+\kappa_{2}$. (Recall, we exclude singular solutions.)
\end{enumerate}
\end{theorem}

\textbf{Proof:} According to Theorem \ref{theorem6} and Corollary \ref{corollary5}, the following possibilities exist:

\begin{enumerate}[(a)]
\item $\tilde{\mathcal{K}}=\bar{L}\mathrm{\thinspace }=L_{1}$.
\item $\tilde{\mathcal{K}}=L_{1}\left[ + \right]\thinspace L_{2}$, where $L_{2}$ is a non-positive (non-degenerate) subspace.
\item $\tilde{\mathcal{K}}=L_{1}\left[ + \right]\thinspace L_{2}$, where $L_{2}$ is a positive subspace.
\item $\tilde{\mathcal{K}}=L_{1}\left[ + \right]\left( L_{0}\dot{+}F \right)\left[ + \right]L_{2}$.
\end{enumerate}
	
(i) By defintion the existence of the nontrivial (which is also non-singular) solution $(z_{1}, z_{2}; h_{1}, h_{2})$ of (\ref{eq52}) means that for at least one function $Q_{i}, \thinspace i=1, \thinspace 2$, it holds: 
\[
h_{i} \ne 0 \wedge \Gamma_{iz_{i}}h_{i} \ne 0. 
\]
In other words, the existence of the nontrivial solution $(z_{1}, z_{2}; h_{1}, h_{2})$ of (\ref{eq52}) is equivalent to: 
\[
\left( {\begin{array}{*{20}c}
\Gamma_{1z_{1}}h_{1}\\
\Gamma_{2z_{2}}h_{2}\\
\end{array} } \right)
\ne 
\left( {\begin{array}{*{20}c}
0\\
0\\
\end{array} } \right).
\]
and
\[
\left( \frac{Q_{1}\left( z_{1} \right)-Q_{1}\left( \bar{w} \right)}{z_{1}-\bar{w}}h_{1},g \right)+\thinspace \left(\frac{Q_{2}\left(z_{2} \right)-Q_{2}\left( \bar{w} \right)}{z_{2}-\bar{w}}h_{2},g \right)
\]
\[
=\left[ \left( {\begin{array}{*{20}c}
\Gamma_{1z_{1}}h_{1}\\
\Gamma_{2z_{2}}h_{2}\\
\end{array} } \right),\left( {\begin{array}{*{20}c}
\Gamma_{1w}g\\
\Gamma_{2w}g\\
\end{array} } \right) \right]=0,\thinspace \forall w\in \mathcal{D}\left( Q 
\right)\thinspace ,\forall g\in \mathcal{H}. 
\]
This is equivalent to existence of 
\[
0\ne \left( {\begin{array}{*{20}c}
\Gamma_{1z_{1}}h_{1}\\
\Gamma_{2z_{2}}h_{2}\\
\end{array} } \right)\in L^{\left[ \bot \right]}.
\]
This is further equivalent to the claim that one of the cases (b), (c), or (d) is satisfied, which is  according to Corollary \ref{corollary5} (i) equivalent to the claim that $\tilde{\mathcal{K}}$ is not minimal state space of $Q$. 

(ii) In Section 3.2 we showed that we can identify $\varepsilon_{z}h=\tilde{\Gamma }_{z}$. Solution $\left( z;h \right)$ is a nontrivial solution of (\ref{eq53}), if and only if $\tilde{\Gamma }_{z}h\in L$ and $\left[ \tilde{\Gamma }_{z}h,\tilde{\Gamma }_{w}g \right]=0, \forall w\in \mathcal{D}\left( Q \right),\\
 \forall g\in \mathcal{H}$. This is equivalent to $0 \neq \tilde{\Gamma }_{z}h \in L \cap L^{[\bot]}$, i.e. it is an isotropic element in $L$.  

(iii) If nontrivial and non-positive solutions of (\ref{eq52}) exist, then (b) or (d) holds. Therefore, $\kappa =\tilde{\kappa }_{1}<\tilde{\kappa }= \kappa_{1}+\kappa_{2}$. 

(iv) Let us first prove the claim: If $G=l.s.\lbrace x:[x,x]>0 \rbrace$, then $G$ is a positive manifold. 
  
If $x$ and $y$ are two positive and linearly dependent vectors, i.e. $y = \beta x, \beta \ne -1$, then obviously $[x +y, x +y]>0$. 

Assume now that $x$ and $y$ are two positive and linearly independent vectors. For every $\alpha = \vert \alpha \vert e^{i \varphi } \in C$ it holds $\vert \alpha \vert ^{2}[x,x]=[\alpha x,\alpha x]=\vert \alpha \vert ^{2} [e^{i \varphi } x,e^{i \varphi} x]>0$. Because of this property, in the sequel we can consider $\alpha \in R$ in the linear combinations of the form $\alpha x + y$, without loss of generality. 

Then, for every $\alpha \in R$ and two positive independent vectors $x,y \in G$, it holds 
\[
P(\alpha):=\left[ \alpha x + y, \alpha x + y \right]=\alpha ^{2}[x,x]+2Re[x,y]\alpha+[y,y]\geq 0.
\]
Quadratic polynomial $P(\alpha)\geq 0$ because $[x,x]>0$ and its discriminant is non-positive, according to Cauchy-Schwartz inequality. As we know, equality sign in Cauchy-Schwartz inequality holds only when  $x$ and $y$ are two linearly dependent vectors. Hence, for  positive independent vectors $x,y \in G$, which we have here, it holds $\left[ \alpha x + y, \alpha x + y \right]>0$. Because we already proved that linear combination of two linearly dependent positive vectors $x$ and $y$ is positive, we can claim that linear combination of any two positive vectors is a positive vector.   

Then the positivity of a linear combination of $n$ positive vectors follows by induction.  
 
Assume now that all solutions of (\ref{eq52}) are positive. According to the above claim, quadratic form in the second sum of (\ref{eq50}) is positive. Then the scalar product in the subspace 
\[
L^{\left[ \bot \right]}=c.l.s.\left\{ \left( {\begin{array}{*{20}c}
\Gamma_{1z_{1}}h_{1}\\
\Gamma_{2z_{2}}h_{2}\\
\end{array} } \right): z_{i} \in \mathcal{D}(Q_{i}), \thinspace \left[ \left( {\begin{array}{*{20}c}
\Gamma_{1z_{1}}h_{1}\\
\Gamma_{2z_{2}}h_{2}\\
\end{array} } \right),\tilde{\Gamma }_{w}g \right]=0,\thinspace \forall w\in 
\mathcal{D}\left( Q \right), \forall g\in \mathcal{H} \right\},
\]
is non-negative or positive definite. We know that $L^{\left[ \bot \right]} = L_{0}[+]L_{2}$, where $L_{0}$ is isotropic subspace of $L^{\left[ \bot \right]}$ and $L_{2}$ is a positive definite subspace, see \cite[Theorem 3.3]{IKL}. Therefore, if there exist a neutral vector $e \in L^{\left[ \bot \right]}$, then it has to be in $L_{0}$. That is equivalent to $\tilde{\kappa_{1}}<\kappa_{1}+\kappa_{2}$. $L_{0}\neq \lbrace 0 \rbrace $ means that $\bar{L}$ is degenerate. According to Lemma \ref{lemma10}, then $L(Q)$ is also degenerate. 

If  $L_{0}=\lbrace 0 \rbrace$ we have case (c), which is equivalent to $\kappa = \kappa_{1}+\kappa_{2}$. 
In Example \ref{example4} we will prove existence of the case (c).

(v) If we assume, in contrast to the first claim of (v), that $(z;h)$ is a nontrivial solution of (\ref{eq53}), then we get 
\[
\left[ \left( {\begin{array}{*{20}c}
\Gamma_{1z}h\\
\Gamma_{2z}h\\
\end{array} } \right),\left( {\begin{array}{*{20}c}
\Gamma_{1w}g\\
\Gamma_{2w}g\\
\end{array} } \right) \right]=0,\thinspace \forall w\in \mathcal{D}\left( Q \right)\thinspace ,\forall g\in \mathcal{H}.
\]
This means that $0\neq \tilde{\Gamma }_{z}h\in L_{0}^{0}$. According to Corollary \ref{corollary5} (iii), it holds $\kappa <\thinspace \kappa_{1}+\kappa_{2}$. This is a contradiction that proves the first claim of (v). 

In Example \ref{example6} we will see that even when equation (\ref{eq52}) has only trivial solution it is possible to have $\kappa < \kappa_{1}+\kappa_{2}$. That will prove the second claim of (v). $\square $
\\

The following theorem gives us some analytic tools to research existence of positive, negative, isotropic and neutral vectors in $\bar{L}$.

\begin{theorem}\label{theorem10} Assume $Q_{i}{\in N}_{\kappa_{i}}\left( \mathcal{H} \right)$ are minimally represented by triplets $\left( \mathcal{K}_{i},A_{i},\Gamma_{i}\thinspace \right), \\
i=1,\thinspace 2,$\textit{ and }$Q:=Q_{1}+Q_{2}$ is represented by $( \tilde{\mathcal{K}},\tilde{A,}\thinspace \tilde{\Gamma } )$.

\begin{enumerate}[(i)]
\item There exist 
\end{enumerate}
\begin{equation}
\label{eq54}
e:=\lim\limits_{n\to\infty}{\tilde{\Gamma }_{z_{n}}h_{z_{n}}}\left( \in \bar{L} \right),\thinspace 
\thinspace e\ne 0,
\end{equation}
if and only if it holds:

\begin{enumerate}[(a)]
\item $\left( \frac{Q_{1}\left( z_{n} \right)-Q_{1}\left( \bar{w_{1}} \right)}{z_{n}-\bar{w_{1}}}h_{z_{n}},g^{1} \right)+\left( \frac{Q_{2}\left( z_{n} \right)-Q_{2}\left( \bar{w_{2}} \right)}{z_{n}-\bar{w_{2}}}h_{z_{n}},g^{2} \right)\to a\left( w_{1},w_{2},g^{1},g^{2} \right)\not\equiv 0,\thinspace \left( n\to \infty \right),\\
\forall w_{i}\in \mathcal{D}\left( Q_{i} \right)\thinspace , \forall g^{i}\in \mathcal{H},\thinspace i=1,\thinspace 2,$ 
\item $\left( \left( \frac{Q_{1}\left( z_{n} \right)-Q_{1}\left( \bar{z_{n}} \right)}{z_{n}-\bar{z_{n}}}+\thinspace \frac{Q_{2}\left( z_{n} \right)-Q_{2}\left( \bar{z_{n}} \right)}{z_{n}-\bar{z_{n}}} \right)h_{z_{n}},h_{z_{n}} \right)\to b\ne \mp \infty ,\thinspace \left( n\to \infty \right)$
\end{enumerate}
for some sequences $\lbrace z_{n} \rbrace _{n=1}^{\infty}\subseteq \mathcal{D}(Q)$, $\lbrace h_{n} \rbrace _{n=1}^{\infty}\subseteq \mathcal{H}$. 

In that case $ e \in \bar{L}$ is positive, neutral, negative element if and only if $ b>,\thinspace =,\thinspace <\thinspace 0$, respectively.

\begin{enumerate}[(ii)]
\item If in addition to (a) and (b) it holds 
\end{enumerate}
\begin{enumerate}[(c)]
\item $\left( \left( \frac{Q_{1}\left( z_{n} \right)-Q_{1}\left( \bar{w} \right)}{z_{n}-\bar{w}}+\frac{Q_{2}\left( z_{n} \right)-Q_{2}\left( \bar{w} \right)}{z_{n}-\bar{w}} \right)h_{z_{n}}, g \right) \to 0,\thinspace \left( n\to \infty \right),\forall w \in \mathcal{D} \left( Q \right), \thinspace \forall g \in \mathcal{H} $,
\end{enumerate}

then element $e$ given by (\ref{eq54}) is an isotropic vector of $\bar{L}, b=0, \tilde{\mathcal{K}}$ is not minimal state space of $Q$, and $\kappa <\thinspace \kappa_{1}+\kappa_{2}$.
\end{theorem}

\textbf{Proof:} (i) Let us assume that (\ref{eq54}) holds. According to \cite[Theorem 2.4]{IKL}, it is equivalent to 
\[
\left[ \tilde{\Gamma }_{z_{n}}h_{z_{n}},\left( {\begin{array}{*{20}c}
\Gamma_{1w_{1}}g^{1}\\
\Gamma_{2w_{2}}g^{2}\\
\end{array} } \right) \right]\to \left[ e,\left( {\begin{array}{*{20}c}
\Gamma_{1w_{1}}g^{1}\\
\Gamma_{2w_{2}}g^{2}\\
\end{array} } \right) \right]\not\equiv 0,\left( n\to \infty 
\right);w_{i}\in \mathcal{D}\left( Q_{i} \right)\thinspace ,g^{i}\in \mathcal{H},\thinspace 
i=1,\thinspace 2,
\]
and
\[
\left[ \tilde{\Gamma }_{z_{n}}h_{z_{n}},\tilde{\Gamma }_{z_{n}}h_{z_{n}} 
\right]\to \left[ e,e \right]\left( n\to \infty \right).
\]
Those limits can be written as
\[
\left( \frac{Q_{1}\left( z_{n} \right)-Q_{1}\left( \bar{w_{1}} 
\right)}{z_{n}-\bar{w_{1}}}h_{z_{n}},g^{1} \right)+\left( \frac{Q_{2}\left( 
z_{n} \right)-Q_{2}\left( \bar{w_{2}} 
\right)}{z_{n}-\bar{w_{2}}}h_{z_{n}},g^{2} \right)
\]
\[
\to \left[ e,\left( {\begin{array}{*{20}c}
\Gamma_{1w_{1}}g^{1}\\
\Gamma_{2w_{2}}g^{2}\\
\end{array} } \right) \right]=:a \not\equiv 0,\left( n\to \infty \right)
\]
and
\[
\left( \left( \frac{Q_{1}\left( z_{n} \right)-Q_{1}\left( \bar{z_{n}} 
\right)}{z_{n}-\bar{z_{n}}}+\thinspace \frac{Q_{2}\left( z_{n} 
\right)-Q_{2}\left( \bar{z_{n}} \right)}{z_{n}-\bar{z_{n}}} 
\right)h_{z_{n}},h_{z_{n}} \right)\to \left[ e,e \right]=:b, \thinspace \left( n\to \infty 
\right).
\]
Because, $b:=[e,e]$, the last statement of (i) holds by definition. This proves (i). 

(ii) Assume now that condition (c) is satisfied as well. Then $\forall w\in \mathcal{D}\left( Q 
\right),\forall g\in \mathcal{H}$ it holds
\[
\left( \left( \frac{Q_{1}\left( z_{n} \right)-Q_{1}\left( \bar{w} 
\right)}{z_{n}-\bar{w}}+\frac{Q_{2}\left( z_{n} \right)-Q_{2}\left( \bar{w} 
\right)}{z_{n}-\bar{w}} \right)h_{z_{n}},g \right)=\left[ \tilde{\Gamma 
}_{z_{n}}h_{z_{n}},\tilde{\Gamma }_{w}g \right]\to \left[ e,\tilde{\Gamma 
}_{w}g \right]=0,
\]
when $\left( n\to \infty \right)$. This means that $e\ne 0$, and $e\in \bar{L}\cap L^{\left[ \bot \right]}=L_{0}$. Hence, it must be $b=0$. According to Corollary \ref{corollary5} (iii) it must be $\kappa <\kappa_{1}+\kappa_{2}$. $\square$
\\

\textbf{4.3} The following simple examples clarify the previous statements. In addition, 
they serve as proofs of existence of the cases theoretically anticipated in 
Corollary \ref{corollary5} and Theorem \ref{theorem8}. 

\begin{example}\label{example4}  Consider the following matrix functions that satisfy 
conditions of Theorem \ref{theorem8}.
\end{example}
\[
Q_{1}\left( z \right)=-\left( {\begin{array}{*{20}c}
z^{-1} & 0\\
0 & z^{-1}\\
\end{array} } \right)\in N_{0}\left( \mathbf{C}^{2} \right),\thinspace 
Q_{2}\left( z \right)=-\left( {\begin{array}{*{20}c}
z^{-2} & z^{-1}\\
z^{-1} & 0\\
\end{array} } \right)\in N_{1}\left( \mathbf{C}^{2} \right).
\]
Then,
\[
Q\left( z \right):
=Q_{1}\left( z \right)+Q_{2}\left( z \right)
=-\left({\begin{array}{*{20}c}
z^{-1}+z^{-2} & z^{-1}\\
z^{-1} & z^{-1}\\
\end{array} } \right).
\]
We can then solve (\ref{eq52}):
\[
\left( {\begin{array}{*{20}c}
\frac{1}{z_{1}\bar{w}} & 0\\
0 & \frac{1}{z_{1}\bar{w}}\\
\end{array} } \right)\left( {\begin{array}{*{20}c}
h_{1}^{1}\\
h_{2}^{1}\\
\end{array} } \right)+\left( {\begin{array}{*{20}c}
\frac{z_{2}+\bar{w}}{z_{2}^{2}\bar{w}^{2}} & \frac{1}{z_{2}\bar{w}}\\
\frac{1}{z_{2}\bar{w}} & 0\\
\end{array} } \right)\left( {\begin{array}{*{20}c}
h_{1}^{2}\\
h_{2}^{2}\\
\end{array} } \right)=0,\thinspace \forall w\in \mathcal{D}\left( Q \right).
\]
Solving this system gives 
\[
h^{1}=\left( {\begin{array}{*{20}c}
h_{1}^{1}\\
0\\
\end{array} } \right),\thinspace h^{2}=\left( {\begin{array}{*{20}c}
0\\
-\thinspace \frac{z_{2}}{z_{1}}h_{1}^{1}\thinspace \\
\end{array} } \right).
\]
Then we easily verify that all nontrivial solutions are positive. Indeed, $E=$
\[
\left( \frac{Q_{1}\left( z_{1} \right)-Q_{1}\left( \bar{z_{1}} 
\right)}{z_{1}-\bar{z_{1}}}\left( {\begin{array}{*{20}c}
h_{1}^{1}\\
0\\
\end{array} } \right),\left( {\begin{array}{*{20}c}
h_{1}^{1}\\
0\\
\end{array} } \right) \right)+\left( \frac{Q_{2}\left( z_{2} 
\right)-Q_{2}\left( \bar{z_{2}} \right)}{z_{2}-\bar{z_{2}}}\left( 
{\begin{array}{*{20}c}
0\\
-\thinspace \frac{z_{2}}{z_{1}}h_{1}^{1}\thinspace \\
\end{array} } \right),\left( {\begin{array}{*{20}c}
0\\
-\thinspace \frac{z_{2}}{z_{1}}h_{1}^{1}\thinspace \\
\end{array} } \right) \right)
\]
\[
=2\frac{\left| h_{1}^{1} \right|^{2}}{\left| z_{1} \right|^{2}}>0.
\]
According to Theorem \ref{theorem8} (i), the representation (\ref{eq36}) is not minimal. 

To determine whether the number of negative squares is preserved we 
can apply Definition \ref{definition2}. We can take $n=1;\thinspace z_{1}\in 
C^{+},\thinspace Re\thinspace z_{1}<0,\thinspace \thinspace h=\left( 
{\begin{array}{*{20}c}
1\\
-1\\
\end{array} } \right)$. Then from 
\[
\left( N_{Q}\left( z_{1},z_{1} \right)h_{1},h_{1} \right)=\frac{2Re(z_{1})}{\left| z_{1} \right|^{4}}<0
\]
and $\kappa_{1}+\kappa_{2}=1$, we conclude $\kappa =1$. Hence, we have 
that all nontrivial solutions are positive, and number of negative squares 
is preserved even though representation (\ref{eq36}) is not minimal. This also proves existence of the case (c) in the proof of Theorem \ref{theorem8}.

We have already proved that $L^{\left[ \bot \right]}$ contains positive elements. Because of $\kappa=\thinspace \kappa_{1}+\kappa_{2}=1$ we know that $\bar{L}$ is non-degenerate. Therefore, $L_{2}:=L^{\left[ \bot \right]}$ is a positive subspace. This proves the existence of the case anticipated in Corollary \ref{corollary5} (ii) and thus completes the proof of Corollary \ref{corollary5} (ii). $\square$

Note that without Theorem \ref{theorem8}, we would have to find operator representations of the functions $Q_{i}$ and $Q$ to obtain the above answers, which would make the task much more difficult. 

\begin{example}\label{example6}  Consider the functions $Q_{1}\left( z 
\right):=-2z^{-1}-z^{-2}\in N_{1}$ and $Q_{2}\left( z \right):=2z^{-1}\in 
N_{1}$. Then
\[
Q\left( z \right):=Q_{1}\left( z \right)+Q_{2}\left( z \right)=-z^{-2}\in 
N_{1}.
\]
Hence, $\kappa_{1}+\kappa_{2}=2>1=\kappa $. 
\end{example}
In this example, (\ref{eq52}) is given by
\[
\left( \frac{2}{z_{1}\bar{w}}+\frac{z_{1}+\bar{w}}{z_{1}^{2}\bar{w}^{2}} 
\right)h_{1}-\frac{2}{z_{2}\bar{w}}h_{2}=0,\thinspace \forall w\in \mathcal{D}\left( Q 
\right).
\]

This equation has only the trivial solution $h_{1}=h_{2}=0$. Hence, this is an example of the sum $Q:=Q_{1}+Q_{2}$ that has only a trivial solution of (\ref{eq52}) and still does not preserve the number of negative squares. This completes the proof of Theorem \ref{theorem8} (v).

According to Theorem \ref{theorem8} the subspace $L^{\left[\bot \right]}$ should be non-positive. We will prove that $L^{\left[ \bot \right]}$ is negative. This will also prove existence of the case (b) in the proof of Theorem \ref{theorem8}. In order to do that we will use operator representations:
\[
Q_{1}\left( z \right):=\Gamma_{1}^{+}\left( A_{1}-zI \right)^{-1}\Gamma_{1}=-2z^{-1}-z^{-2}\in N_{1},
\]
where 
\[
 A_{1}=\left( {\begin{array}{*{20}c}
 0 & 1\\
 0 & 0\\
 \end{array} } \right),\thinspace \thinspace J_{1}=\left( 
 {\begin{array}{*{20}c}
 0 & 1\\
 1 & 0\\
 \end{array} } \right),\thinspace \Gamma_{1}
=\left( {\begin{array}{*{20}c}
1\\
1\\
\end{array} } \right),\thinspace \Gamma_{1}^{+}
=\Gamma_{1}^{\ast}J_{1}=\left( {\begin{array}{*{20}c}
1 & 1\\
\end{array} } \right)\left( {\begin{array}{*{20}c}
0 & 1\\
1 & 0\\
\end{array} } \right)
\]
 and
\[
Q_{2}\left( z \right):=\Gamma_{2}^{+}\left( A_{2}-zI \right)^{-1}\Gamma_{2}=2z^{-1}\in N_{1}
\]
where
\[
A_{2}=\left( 0 \right),\thinspace J_{2}=\left( -1 \right),\thinspace \Gamma 
_{2}=\left( \sqrt[2]{2} \right),\thinspace \Gamma_{2}^{+}=\Gamma_{2}^{\ast 
}J_{2}=-\left( \sqrt[2]{2} \right).
\]
According to the definitions in Section \ref{s6} we have,
\[
\tilde{A}=\left( {\begin{array}{*{20}c}
0 & 1 & 0\\
0 & 0 & 0\\
0 & 0 & 0\\
\end{array} } \right),\thinspace \tilde{J}
=\left( {\begin{array}{*{20}c}
0 & 1 & 0\\
1 & 0 & 0\\
0 & 0 & -1\\
\end{array} } \right),\thinspace \tilde{\Gamma }
=\left( {\begin{array}{*{20}c}
1\\
1\\
\sqrt 2 \\
\end{array} } \right),
\]
\[
\tilde{\Gamma }^{+}
=\left( {\begin{array}{*{20}c}
1 & 1 & \sqrt 2 \\
\end{array} } \right)\left( {\begin{array}{*{20}c}
0 & 1 & 0\\
1 & 0 & 0\\
0 & 0 & -1\\
\end{array} } \right),
\]
\[
\tilde{\Gamma }_{z}={\thinspace \left( \tilde{A}-zI 
\right)}^{-1}\tilde{\Gamma }h=-\left( {\begin{array}{*{20}c}
z^{-1} & z^{-2} & 0\\
0 & z^{-1} & 0\\
0 & 0 & z^{-1}\\
\end{array} } \right)\left( {\begin{array}{*{20}c}
1\\
1\\
\sqrt 2 \\
\end{array} } \right)h=-\left( {\begin{array}{*{20}c}
z^{-1}+z^{-2}\\
z^{-1}\\
\sqrt 2 z^{-1}\\
\end{array} } \right)h\in L,
\]
and 
\[
\bar{L}=c.l.s. \left\{ -\left( {\begin{array}{*{20}c}
z^{-1}+z^{-2}\\
z^{-1}\\
\sqrt 2 z^{-1}\\
\end{array} } \right)h:\thinspace z\in \mathcal{D}\left( Q \right),\thinspace h\in \mathcal{H} \right\}.
\]
Then for $y=\left( {\begin{array}{*{20}c}
y_{1}\\
y_{2}\\
y_{3}\\
\end{array} } \right)\in {L^{[\bot ]}\mathrm{\thinspace \subseteq 
}\tilde{\mathcal{K}}=\mathcal{K}}_{1}\left[ + \right]\mathcal{K}_{2}$ we have
\[
0=\left[ y,\left( \tilde{A}-zI \right)^{-1}\tilde{\Gamma }h \right]=\left( 
\left( {\begin{array}{*{20}c}
y_{1}\\
y_{2}\\
y_{3}\\
\end{array} } \right),\tilde{J}\left( {\begin{array}{*{20}c}
{-z}^{-1}-z^{-2}\\
-z^{-1}\\
-\sqrt 2 z^{-1}\\
\end{array} } \right)h \right)
\]
\[
=\left( \left( {\begin{array}{*{20}c}
y_{1}\\
y_{2}\\
y_{3}\\
\end{array} } \right),\left( {\begin{array}{*{20}c}
{-z}^{-1}\\
{-z}^{-1}-z^{-2}\\
\sqrt 2 z^{-1}\\
\end{array} } \right)h \right),\thinspace \forall z\in \mathcal{D}\left( Q 
\right),\forall h\in \mathcal{H}.
\]
\[
\Rightarrow y=\left( {\begin{array}{*{20}c}
y_{1}\\
0\\
y_{1} \mathord{\left/ {\vphantom {y_{1} \sqrt 2 }} \right. 
\kern-\nulldelimiterspace} \sqrt 2 \\
\end{array} } \right).
\]
\[
\left[ y,y \right]=\left( \left( {\begin{array}{*{20}c}
y_{1}\\
0\\
y_{1} \mathord{\left/ {\vphantom {y_{1} \sqrt 2 }} \right. 
\kern-\nulldelimiterspace} \sqrt 2 \\
\end{array} } \right),\tilde{J}\left( {\begin{array}{*{20}c}
y_{1}\\
0\\
y_{1} \mathord{\left/ {\vphantom {y_{1} \sqrt 2 }} \right. 
\kern-\nulldelimiterspace} \sqrt 2 \\
\end{array} } \right) \right)=-\frac{\left| y_{1} \right|^{2}}{2}<0.
\]
Hence, vector $\mathbf{y}\in L^{\left[ \bot \right]}$ is strictly negative. This is indeed the case (b) anticipated in the proof of Theorem \ref{theorem8}. $\square$

\section{ The final decomposition of Q}\label{s10}

\textbf{5.1}  Let the function $Q\in N_{\kappa }(\mathcal{H})$ be minimally represented by (\ref{eq2}) and let $\alpha \in R$ be a generalized pole of $Q$ that is not of positive type. It is customary to say that $A$ and $\Gamma $ are \textit{closely connected} if representation (\ref{eq2}) is minimal. Let us decompose the function $Q$ by means of the Jordan chains of the representing relation $A$ at $\alpha $. 

According to \cite[Lemma 1]{B}, there is no loss of generality to assume that $\alpha \in R$ is a single generalized pole that is not of positive type. In that case $A$ is an operator. For given eigenvector $x_{0}$ of $A$ at $\alpha \in R$, let us denote by $X$ one of the maximal Jordan chains of $x_{0}$. Let us denote by 
\[
S_{\alpha }\left( x_{0} \right):=l.s.\left\{ X \right\}.
\]
Let the Hilbert subspace, denoted here by $\mathcal{K}_{0} \subset \mathcal{K}$, consist of all positive 
eigenvectors of the representing operator $A$ at $\alpha $. Let $E_{0}:\mathcal{K}\to 
\mathcal{K}_{0}$ be the orthogonal projection $E^{'}:=I-E_{0}$, $\mathcal{K}^{'}:=E^{'}\mathcal{K}$ and 
$\Gamma_{0}:=E_{0}\Gamma $. Subspaces $\mathcal{K}_{0}$ and $\mathcal{K}^{'}$ obviously 
reduce operator $A$. We define $\Gamma^{'}:=E^{'}\Gamma $ and 
$A^{'}:=E^{'}AE^{'}$. 

Now let $x_{0}^{1},\thinspace \mathellipsis ,\thinspace x_{l_{1}-1}^{1}$ be a 
maximal non-degenerate Jordan chain of $A^{'}$ at $\alpha $ in the 
Pontryagin space $\mathcal{K}^{'}$. We define the projection: $E_{1}:\mathcal{K}^{'}\to 
S_{\alpha }(x_{0}^{1})$, and sub-space ${\mathcal{K}_{1}:=E}_{1}\mathcal{K}^{'}$. Then 
$A_{1}=E_{1}A^{'}E_{1}$ and $\Gamma_{1}:=E_{1}\mathrm{\thinspace }\Gamma 
^{'}$ are closely connected operators. Let $\thinspace \kappa_{1}$ denote the 
negative index of the Pontryagin space $\mathcal{K}_{1}$. 

We can repeat these steps until we exhaust all non-degenerate Jordan chains. 
At every step we can decompose the corresponding function as in Theorem \ref{theorem3}. 

Assume that there are $r>0$ such (non-degenerate) chains at $\alpha $. We 
introduce $E:=E_{0}+E_{1}+\mathellipsis +E_{r}$. Then, $\mathcal{K}=E\mathcal{K}\left[ + 
\right]\left( I-E \right)\mathcal{K}$. Let us introduce $E_{r+1}:=I-E$, 
$\mathcal{K}_{r+1}:=E_{r+1}\mathcal{K}$, $\Gamma_{r+1}=E_{r+1}\Gamma $. Subspaces $E\mathcal{K}$ and 
$\mathcal{K}_{r+1}$ obviously reduce $A$. From the construction of the Pontryagin 
space $\mathcal{K}_{r+1}$ we conclude that all degenerate chains of $A$ at $\alpha $ 
are in $\mathcal{K}_{r+1}$. 

By using the above notation, we can summarize these results in the following 
proposition.

\begin{proposition}\label{proposition6} Let $\alpha \in R$ be a generalized pole that is not of positive type of $Q\in N_{\kappa }(\mathcal{H})$, where $Q$ is given by minimal representation (\ref{eq2}). Then 
\begin{equation}
\label{eq62}
\ \mathcal{K}=\mathcal{K}_{0}\left[ + \right]\mathcal{K}_{1}\left[ + \right]\mathellipsis \left[ 
+ \right]\mathcal{K}_{r}\left[ + \right]\mathcal{K}_{r+1}
\end{equation}
where $r\in \mathbf{N}$ is the number of independent non-degenerate Jordan chains of $A$ at $\alpha ; \mathcal{K}_{i}$ are A-invariant Pontryagin subspaces of indices $\kappa_{i}, \thinspace i=0,\thinspace 1,\thinspace \mathellipsis ,\thinspace r,\thinspace r+1$, respectively; $\thinspace \kappa 
_{0}=0, \kappa =\sum\limits_{i=1}^{r+1} \kappa_{i}$. For every $i=1,\thinspace 
2,\thinspace \mathellipsis ,\thinspace r$, subspace $\mathcal{K}_{i}$ is a linear span of the corresponding maximal non-degenerate Jordan chain $\thinspace 
x_{0}^{i},\thinspace \mathellipsis ,\thinspace x_{l_{i}-1}^{i}$. All positive eigenvectors are in $\mathcal{K}_{0}$. 
All degenerate chains of $A$ at $\mathrm{\alpha }$ are in $\mathcal{K}_{r+1}$. 

The corresponding nontrivial decomposition $Q:=Q_{0}+Q_{1}+\mathellipsis +Q_{r}+Q_{r+1}$ satisfies $\kappa 
=\sum\limits_{i=0}^{r+1} \kappa_{i}.$
\end{proposition}

\textbf{5.2} Because $\mathcal{K}_{i}$, $i=1,\thinspace \mathellipsis ,\thinspace r$, is a linear span of a maximal Jordan chain, it does not have a nontrivial invariant subspaces of $A$. In Proposition \ref{proposition6}, we separated non-degenerate maximal Jordan chains $X^{i}$ and $X^{j}$ by $A$-invariant disjoint subspaces $\mathcal{K}_{i}$ and $\mathcal{K}_{j}$, i.e. $X^{i} \subset \mathcal{K}_{i}, \thinspace X^{j} \subset \mathcal{K}_{j}$, $\mathcal{K}_{i} \cap \mathcal{K}_{j}=\left\{ 0 \right\},\thinspace \forall i\ne j$. The following natural question arises: Is it possible to separate \textbf{degenerate} Jordan chains in a similar way? More precisely:

Let $A$ be a self-adjoint operator in a Pontryagin space $\mathcal{K}$. \textit{Given two degenerate maximal Jordan chains} $X^{i},\thinspace i=1,\thinspace 2,$ \textit{at an eigenvalue} $\alpha \in R$, \textit{is it possible to find an} $A$\textit{-invariant non-degenerate subspace} $\mathcal{K}_{1}$ \textit{such that it holds} $X^{1}\subset \mathcal{K}_{1}$ \textit{and} $\thinspace X^{2}\cap \mathcal{K}_{1}=\emptyset ?$

In order to address this question, we introduce the following model with two 
independent degenerate chains at $\alpha =0$ of the first order, i.e. two 
neutral eigenvectors. We denote: 
$\langle k \rangle=l.s.\lbrace k \rbrace$.

\begin{proposition}\label{proposition8} Assume that

\begin{equation}
\label{eq64}
\mathcal{K}=\mathcal{H}\left[+ \right]\left(\left(\langle x_{0}^1 \rangle \left[+ \right]\langle x_{0}^2 \rangle \right)\dotplus\left(\langle f^1 \rangle \left[+ \right]\langle f^2 \rangle \right)\right),
\end{equation}

\[
J=\left({\begin{array}{l}
{\begin{array}{*{20}c}
I & 0 & 0 & 0 & 0\\
0 & 0 & 0 & 1 & 0\\
0 & 0 & 0 & 0 & 1\\
0 & 1 & 0 & 0 & 0\\
0 & 0 & 1 & 0 & 0\\
\end{array}} \\
\end{array}} \right),\thinspace \thinspace A=\left({\begin{array}{l}
{\begin{array}{*{20}c}
A_{11} & 0 & 0 & a_{1} & a_{2}\\
\left(.,a_{1} \right) & 0 & 0 & \alpha_{1} & 0\\
\left(.,a_{2} \right) & 0 & 0 & 0 & \alpha_{2}\\
0 & 0 & 0 & 0 & 0\\
0 & 0 & 0 & 0 & 0\\
\end{array}}
\end{array}} \right),
\]
where $\left( \mathcal{H},\left( .,. \right) \right)$ is a Hilbert space, $A_{11}$ is 
a bounded self-adjoint operator on $\mathcal{H}$, ${0\ne \alpha }_{i}\in R,\thinspace 
0\mathbf{\ne }a_{i}\in \mathcal{H},\thinspace i=1,\thinspace 2$, are linearly 
independent. Then 

\begin{enumerate}[(i)]
\item Operator $A$ is a self-adjoint operator in the Pontryagin space $\mathcal{K}$. 
\item Vectors $x_{0}^{i}$ are neutral, simple eigenvectors of $A$ at $\alpha =0$ and $f^{i}=Jx_{0}^{i},\thinspace i=1,\thinspace 2$. 
\item If operator $A_{11}:\mathcal{H}\to \mathcal{H}$ is irreducible, then operator $A$ does not have any eigenvalues different from $\alpha=0$.
\item If operator $A_{11}$ is irreducible, then operator $A$ does not have any invariant non-degenerate subspace that contains one eigenvector $x_{0}^{i}$ and not the other, $x_{0}^{j}, i \neq j; \thinspace i,j = 1, 2$. 
\end{enumerate}
\end{proposition}

\textbf{Proof:} For vectors from $\mathcal{K}$ we will use notation $\left( h,\beta_{1},\beta_{2},\gamma_{1},\gamma_{2}\thinspace \thinspace \right)^{T}, h \in \mathcal{H}, \gamma_{i}, \beta_{i} \in C, i=1, 2$. 
\

Statements (i) and (ii) are straightforward verification. 

(iii) In contrast to the statement, assume that the operator $A$ has the eigenvalue $\beta \ne 0$. Then 
\[
A\left({\begin{array}{l}
{\begin{array}{*{20}c}
h\\
\beta_{1}\\
\beta_{2}\\
\gamma_{1}\\
\gamma_{2}\\
\end{array}}\end{array}}\right)= \beta \left({\begin{array}{l}
{\begin{array}{*{20}c}
h\\
\beta_{1}\\
\beta_{2}\\
\gamma_{1}\\
\gamma_{2}\\
\end{array}}\end{array}}\right)\Rightarrow \left({\begin{array}{l}
{\begin{array}{*{20}c}
A_{11}h+a_{1} \gamma_{1}+a_{2} \gamma_{2}\\
\left(h,a_{1}\right)+\alpha_{1} \gamma_{1}\\
\left(h,a_{2} \right)+\alpha_{2} \gamma_{2}\\
0\\
0\\
\end{array}} \\
\end{array}} \right)=\beta \left({\begin{array}{l}
{\begin{array}{*{20}c}
h\\
\beta_{1}\\
\beta_{2}\\
\gamma_{1}\\
\gamma_{2}\\
\end{array}}\end{array}}\right).
\]
Hence, $\gamma_{i}=0, i=1, 2$. Therefore,

\[
A \left({\begin{array}{l}
{\begin{array}{*{20}c}
h\\
\beta_{1}\\
\beta_{2}\\
0\\
0\\
\end{array}}\end{array}}\right)= \left({\begin{array}{l}
{\begin{array}{*{20}c}
A_{11}h\\
\left(h,a_{1}\right)\\
\left(h,a_{2} \right)\\
0\\
0\\
\end{array}} \\
\end{array}} \right)=\beta \left({\begin{array}{l}
{\begin{array}{*{20}c}
h\\
\beta_{1}\\
\beta_{2}\\
0\\
0\\
\end{array}}\end{array}}\right)\Rightarrow h \ne 0 \thinspace and  \thinspace A_{11}h= \beta h.
\]
This means that operator $A_{11}$ has a nonzero eigenvalue $\beta $. Because 
$A_{11}$ is a bounded self-adjoint operator on $\mathcal{H}$, it is reduced by the eigenvector 
$h\in \mathcal{H}$, see also definition of reductibility in \cite[Section 40]{AG}. That contradicts the assumption that $A_{11}$ is irreducible. This proves (iii). 

(iv) In contrast to the statement, assume that operator $A_{11}$ is 
irreducible in $\mathcal{H}$ and that there exist a non-degenerate, $A$-invariant, 
nontrivial subspace $\mathcal{K}_{1}$ of $\mathcal{K}$ such that $x_{0}^{1}\in \mathcal{K}_{1}$ and 
$x_{0}^{2}\notin \mathcal{K}_{1}$. Then $\mathcal{K}_{1}$ must contain $f^{1}$; otherwise, 
according to (\ref{eq64}), the subspace $\mathcal{K}_{1}$ would be degenerate. Similarly, $\mathcal{K}_{1}$ cannot contain $f^{2}$, because then $\mathcal{K}_{1}$ without $x_{0}^{2}$ would be degenerate. Hence, vectors from $\mathcal{K}_{1}$ must satisfy $\gamma_{1}\neq 0$,$\gamma_{2}= 0$ , and $\mathcal{H}_{1}:=\mathcal{H}\cap \mathcal{K}_{1}$ must contain vectors of the form $A_{11}h+a_{1}\gamma_{1}\in \mathcal{H}$. 
This means, $\mathcal{K}_{1}$ is of the form:
\begin{equation}
\label{eq66}
\mathcal{K}_{1}=\mathcal{H}_{1}\left[+ \right] \left(\langle x_{0}^{1} \rangle \dotplus \langle f^{1} \rangle \right),
\end{equation}
where $\mathcal{H}_{1}\ne \left\{ 0 \right\}$. It is easy to verify, that it holds $x_{0}^{2}\in 
\mathcal{K}_{1}^{\left[ \bot \right]}$. 

For an arbitrarily selected $k_{1}\in \mathcal{K}_{1}$ we have:
\[
k_{1}:=\left({\begin{array}{l}
{\begin{array}{*{20}c}
h\\
\beta_{1}\\
0\\
\gamma_{1}\\
0\\
\end{array}}\end{array}}\right) \in \mathcal{K}_{1}
\Rightarrow A\left({\begin{array}{l}
{\begin{array}{*{20}c}
h\\
\beta_{1}\\
0\\
\gamma_{1}\\
0\\
\end{array}}\end{array}}\right)=\left({\begin{array}{l}
{\begin{array}{*{20}c}
A_{11}h+a_{1} \gamma_{1}\\
\left(h,a_{1}\right)+\alpha_{1} \gamma_{1}\\
\left(h,a_{2} \right)\\
0\\
0\\
\end{array}} \\
\end{array}} \right) \in \mathcal{K}_{1}; \thinspace
 \beta_1, \gamma_{1} \in C, h \in \mathcal{H}_{1}.
\]
Because $\mathcal{K}_{1}$ is $A$-invariant, $Ak_{1}$ must be of the form (\ref{eq66}). Hence, 
it must be
\[
\left( h,a_{2} \right)=0,\forall h\in \mathcal{H}_{1}
\]
and
\[
\thinspace A_{11}h+a_{1}\gamma_{1}\in \mathcal{H}_{1},\thinspace \forall h\in 
\mathcal{H}_{1},\thinspace \forall \gamma_{1}\in \mathbf{C\thinspace }\thinspace .
\]
Hence $a_{2}\bot \mathcal{H}_{1}$, where $0\ne a_{2}\in \mathcal{H}$. This means $\lbrace 0 \rbrace \subsetneq \mathcal{H}_{1} \subsetneq \mathcal{H}$ i.e. $\mathcal{H}_{1}$ is a nontrivial subspace of $\mathcal{H}$. 

If we set $\gamma_{1}=0$ in the second equation, then we conclude that 
$A_{11}h\in \mathcal{H}_{1},\thinspace \forall h\in \mathcal{H}_{1}$. Therefore, $\mathcal{H}_{1}$ is an 
$A_{11}$-invariant, nontrivial subspace in $\mathcal{H}$. Because $A_{11}$ is bounded 
self-adjoint operator on the Hilbert space $\mathcal{H}$, operator $A_{11}$ is reduced 
by $\mathcal{H}_{1}$, see again \cite[Section 40]{AG}. That contradicts the assumption of irreducibility of $A_{11}$ 
and proves (iv).$\square$

This example shows that there does not exist an $A$-invariant subspace that contains one and not the other degenerate eigenvector of $A$ at $\alpha$. 

\begin{corollary}\label{corollary10} There is no nontrivial decomposition of $\mathcal{K}_{r+1}$ and $Q_{r+1}$, i.e. decomposition (\ref{eq62}) of $\mathcal{K}$ and corresponding decomposition of $Q$ are final.
\end{corollary}

\end{document}